\input ./style/arxiv-general.cfg
\documentclass[aap,MSNbibl,seceqn,dvips]{arximspdf}
\makeatletter
   \@ifpackageloaded{graphicx}{}{\usepackage{graphicx}}
\makeatother
\usepackage[mathscr]{euscript}

\doi{10.1214/14-AAP1074}
\volume{25}
\issue{6}
\pubyear{2015}
\firstpage{3295}
\lastpage{3337}
\docsubty{FLA}

\makeatletter
\newcommand{\eqref}[1]{(\ref{#1})}

\newtheorem{theorem}{Theorem}
\newtheorem{lemma} [theorem] {Lemma}
\newtheorem{claim} [theorem] {Claim}
\newtheorem{corollary} [theorem] {Corollary}
\newproclaim{definition} [theorem] {Definition}
\newproclaim{remark} [theorem] {Remark}
\newtheorem{conjecture} [theorem] {Conjecture}
\newtheorem{question} [theorem] {Question}

\newcommand{\eps}{\varepsilon}

\newcommand{\Ccal}[0]{{{\mathcal C}}}
\newcommand{\Dcal}[0]{{{\mathcal D}}}
\newcommand{\Mcal}[0]{{{\mathcal M}}}
\newcommand{\Ascr}[0]{{{\mathscr A}}}
\newcommand{\Cscr}[0]{{{\mathscr C}}}
\newcommand{\Dscr}[0]{{{\mathscr D}}}
\newcommand{\Kscr}[0]{{{\mathscr K}}}
\newcommand{\Mscr}[0]{{{\mathscr M}}}
\newcommand{\Nscr}[0]{{{\mathscr N}}}

\newcommand{\eR}[0]{{ \mathbb R}}
\newcommand{\eN}[0]{{ \mathbb N}}
\newcommand{\Zed}[0]{{ \mathbb Z}}
\newcommand{\Pee}[0]{{{\mathbb P}}}
\newcommand{\Ee}[0]{{{\mathbb E}}}
\newcommand{\isd}[0]{\stackrel{\mathrm{d}}{=}}
\newcommand{\Gtil}[0]{{\tilde{G}}}
\newcommand{\area}{\operatorname{area}}
\newcommand{\Bi}{\operatorname{Bi}}
\newcommand{\Po}{\operatorname{Po}}
\newcommand{\geom}{\operatorname{Geom}}
\newcommand{\Var}{\operatorname{Var}}
\newcommand{\dd}{\,\mathrm{d}}

\newcommand{\rtil}{{\tilde{r}}}

\newcommand{\strat}[1]{{(STR-#1)}}
\newcommand{\Moat}[0]{ {\Mscr}}
\newcommand{\Nbad}[0]{ {N_{\mathrm{bad}}}}
\newcommand{\Ndang}[0]{ {N_{\mathrm{dang}}}}
\newcommand{\Rpsd}[0]{ {R_{\mathrm{psd}}}}
\newcommand{\pdang}[0]{ {p_{\mathrm{dang}}}}
\newcommand{\King}[0]{ {\Kscr}}
\newcommand{\Occ}{\mathrm{Occ}}
\newcommand{\super}{\mathrm{Big}}

\newcommand{\cC}{{\mathcal C}}
\newcommand{\cX}{{\mathcal X}}

\newcommand{\kmin}{{m_{\min}}}
\renewcommand{\rho}{\varrho}
\newcommand{\aas}{a.a.s. }
\newcommand{\seq}{\subseteq}

\newcommand{\geBy}[1]{\stackrel{\scriptsize{\eqref{#1}}}{\ge}}
\newcommand{\leBy}[1]{\stackrel{\scriptsize{\eqref{#1}}}{\le}}
\newcommand{\eqBy}[1]{\stackrel{\scriptsize{\eqref{#1}}}{=}}
\newcommand{\leByM}[1]{\stackrel{#1}{\le}}

\makeatother

\begin{document}
\begin{frontmatter}

\title{A geometric Achlioptas process}
\runtitle{A geometric Achlioptas process}

\begin{aug}
\author[A]{\fnms{Tobias} \snm{M\"uller}\corref{}\thanksref{T1}\ead[label=e1]{t.muller@uu.nl}}
\and
\author[B]{\fnms{Reto} \snm{Sp\"ohel}\thanksref{T2}\ead[label=e2]{reto.spoehel@bfh.ch}}
\runauthor{T.~M\"uller and R.~Sp\"ohel}
\thankstext{T1}{Supported in part by a VENI grant from
Netherlands Organisation for Scientific Research (NWO).}
\thankstext{T2}{Supported in part by a grant of the Swiss National
Science Foundation
while the author was at MPI Saarbr\"ucken.}

\affiliation{Utrecht University and
Bern University of Applied Sciences}

\address[A]{Deparment of Mathematics\\
Utrecht University\\
P.O.~Box 80010\\
3508TA Utrecht\\
The Netherlands\\
\printead{e1}}

\address[B]{BFH-TI \\
Bern University of Applied Sciences\\
Pestalozzistrasse 20\\
3400 Burgdorf\\
Switzerland\\
\printead{e2}}
\end{aug}
%
%

%
\received{\smonth{4} \syear{2013}}
%
\revised{\smonth{5} \syear{2014}}

%
\begin{abstract}
The \emph{random geometric graph} is obtained by sampling $n$ points
from the unit square (uniformly at random and independently), and
connecting two points whenever their distance is at most $r$, for some
given $r=r(n)$.
We consider the following variation on the random geometric graph: in
each of~$n$ rounds in total, a player is offered \emph{two} random
points from the unit square, and has to select exactly one of these two
points for inclusion in the evolving geometric graph.

We study the problem of avoiding a linear-sized (or ``giant'')\vspace*{1pt}
component in this setting. Specifically, we show that for any
$r \ll\break (n \log\log n)^{-1/3}$ there is a strategy that succeeds in
keeping all component sizes sublinear, with probability tending to one
as $n\to\infty$. We also show that this is tight in the following
sense: for any $r \gg(n \log\log n)^{-1/3}$, the player will be forced
to create a component of size $(1-o(1))n$, no matter how he plays,
again with probability tending to one as $n \to\infty$.
We also prove that the corresponding offline problem exhibits a similar
threshold behaviour at $r(n)=\Theta(n^{-1/3})$.

These findings should be compared to the existing results for the
(ordinary) random geometric graph: there a giant component arises with
high probability once $r$ is of order $n^{-1/2}$. Thus, our results
show, in particular, that in the geometric setting the power of choices
can be exploited to a much larger extent than in the classical Erd\H
{o}s--R\'enyi random graph, where the appearance of a giant component can
only be delayed by a constant factor.
\end{abstract}

%

%
\end{frontmatter}

\section{Introduction}\label{sec1}

The \emph{random geometric graph} with parameters $n$ and $r$
is obtained by sampling $n$ points from the unit square (uniformly at
random and independently), and connecting two points whenever
their distance is at most $r$.
The study of this model essentially goes back to Gilbert \cite
{Gilbert61} who defined a very similar model
in 1961; for this reason it is sometimes also called the \emph{Gilbert
random graph}.

Random geometric graphs form an interesting and rich subject from a
purely theoretical perspective, but they are also studied in
relation to a variety of applications.
They have, for instance, been used to model wireless networks (see,
e.g., \cite{wireless}), the growth of tumors \cite{tumours},
protein--protein interactions \cite{protein}, fiber-based
materials \cite
{materials} and many more phenomena.
Random geometric graphs have been the subject of considerable research
effort over the past decades, and quite precise results are now known
for this model on aspects such as connectivity, Hamilton
cycles, the clique number, the chromatic number and random walks on the graph.
(See, e.g., \cite{penrosekconn,BBKMW,twopoint,McDiarmidMullerRGG,CooperFriezeCoverRgg}.)
A comprehensive overview of the results prior to 2003 can be found in
the monograph \cite{PenroseBoek}.

By the results in Chapter~10 of \cite{PenroseBoek} (which build on the
work of several
previous authors, including Gilbert \cite{Gilbert61}), there is a
constant $\lambda_{\mathrm{crit}}$ such that
if $r = \sqrt{\lambda/n}$ with $\lambda\leq\lambda_{\mathrm{crit}}$
then the largest component contains a sublinear
proportion of all vertices, while if $\lambda> \lambda_{\mathrm
{crit}}$ then the largest component
contains a linear fraction of all vertices.
Phrased differently, a ``giant'' component suddenly emerges when the
average degree exceeds a certain constant ($\pi\lambda_{\mathrm
{crit}}$ to be precise).
An interesting detail is that the precise value of $\lambda_{\mathrm
{crit}}$ remains unknown to this date.

\subsection{Our results}

We consider a \emph{power of choices version} of the random geometric
graph. By this, we mean the following probabilistic process:
There are $n$ \emph{rounds} in total, and in each round a player is
offered \emph{two} random points from the unit square, of which he has
to select exactly one for inclusion in the evolving geometric graph.

The objective of the player is to keep the size of the largest
component as small as possible. (In Section~\ref{sec:concl}, we briefly
discuss the setup when the player wants to maximize
the size of the largest component.) In particular, we are interested in
the question for which functions $r=r(n)$ the player can avoid the
formation of a \emph{linear-sized} (``giant'') component with high
probability.\setcounter{footnote}{2}\footnote{We say that a sequence of events $(A_n)_n$ holds
\emph{with high probability} (abbreviated: w.h.p.) if $\Pee( A_n ) =
1-o(1)$ as $n\to\infty$.
As is common in the random graphs literature, we shall often be a
little bit sloppy with our notation and say things like $X_n = o(f(n))$
w.h.p, where $(X_n)_n$ is some sequence of random variables.
This will of course mean that there exists a function $g(n)=o(f(n))$
such that $X_n \leq g(n)$ w.h.p.
The interpretations of $X_n = (1-o(1))f(n)$ w.h.p., $X_n = \Omega
(f(n))$ w.h.p., etc., are analogous.
} This question is answered by the following theorem.

\begin{theorem}\label{thm:online}
Consider the geometric power of choices process defined above.
There exist functions $f, g \dvtx (0,\infty) \to(0,1)$, $g<f$, such that
the following holds:

 If $r = \sqrt[3]{\frac{c}{n \cdot\log\log n}}$ for some
fixed $c > 0$, then:
\begin{longlist}[(ii)]
\item[(i)] There exists a strategy such that, if the player
follows this strategy, then w.h.p. in round $n$ all components are
smaller than $f(c) n$.
\item[(ii)] W.h.p., no matter which strategy the player
employs, the
largest component in round $n$ has order at least $g(c) n$.
\end{longlist}
Moreover, $f(c) \to0$ as $c\downarrow0$ and
$g(c) \to1$ as $c \to\infty$.
\end{theorem}

Here and in the rest of the paper, $\log(\cdot)$ denotes the \emph{base 2}
logarithm (i.e., the inverse of $2^x$).
The functions $f, g$ provided by the proof satisfy $f(c) = O( \sqrt{c}
)$ as $c\downarrow0$ and
$g(c) = 1 - O( 1/\log c )$ as $c\to\infty$.

Theorem~\ref{thm:online} extends to the setting with an arbitrary fixed
number $d\geq2$ of choices per step; the expression for $r$ then needs
to be replaced by
$r = \sqrt[d+1]{c/(n \cdot(\log\log n)^{d-1})}$. We will come back to
this at the end of the paper.

Note that Theorem~\ref{thm:online} implies in particular that for any
$r \ll(n\log\log n)^{-1/3}$, w.h.p. the largest component is
sublinear in size, and that for any $r \gg (n\log\log n)^{-1/3}$,
w.h.p. the largest component will be of size $(1-o(1))n$ no matter how
the player plays.\vspace*{1pt} Thus, Theorem~\ref{thm:online} establishes a
``threshold'' of $\Theta((n\log\log n)^{-1/3})$ for the appearance of a
giant component in the geometric power of choices process. Note that
this threshold is higher than the threshold for the original geometric
random graph by a power of $n$. This is in stark contrast with the
power of choices version of the well-known Erd\H{o}s--R{\'{e}}nyi
process (see Section~\ref{sec:Background} below); there, the appearance
of a giant component can only be delayed by a constant factor.

Observe also that the threshold behaviour is very different from that
of the original random geometric graph setting: Theorem~\ref
{thm:online} states that qualitatively the behaviour of the process is
the same for all $c>0$.
In the standard geometric graph on the other hand, the size of the
largest component jumps from $\Theta(\log n)$ just below the threshold
$r_{\mathrm{crit}}:= \sqrt{\lambda_{\mathrm{crit}}/n}$ to $\Theta(n)$
just above $r_{\mathrm{crit}}$.
(In the well-known Erd\H os--R\'enyi graph, a similar phenomenon
occurs; see, e.g., \cite{randomgraphs}.)
In the process considered in Theorem~\ref{thm:online}, however, the
order of the largest component under optimal play is
$\Theta(n)$ for every $c>0$.
So, in particular, there is no $c_{\mathrm{crit}}$ analogous to
$\lambda
_{\mathrm{crit}}$.
See Section~\ref{sec:concl} for some additional discussion on the
order of the largest component when $r$ is below the threshold $\Theta(
(n \log\log n)^{-1/3} )$.

\textit{The offline setting.}
In the process we discussed so far, the player does not know which
points will arrive in future rounds. It is interesting to consider what
would happen if the player were \emph{clairvoyant}, that is, if he already
knew from the start of the game where all the points in all the rounds
will fall. Put differently, he is given $n$ pairs of random points all
at once, and needs to select one point from each of these pairs.
This is often called the \emph{offline} version of the game, and the
original version is called the \emph{online} version.

Intuitively, the additional advantage of being clairvoyant should allow
the player to delay the onset of a giant component even further. The
next theorem shows that this is indeed true, but that the advantage is
rather modest---the threshold only increases by a factor of $(\log
\log
n)^{1/3}$.

\begin{theorem} \label{thm:offline}
Consider the geometric offline power of choices setting defined above.
There exist functions $f, g \dvtx (0,\infty) \to(0,1)$, $g<f$, such that
the following holds:

 If $r = \sqrt[3]{\frac{c}{n}}$ for some fixed $c > 0$, then:
\begin{longlist}[(ii)]
\item[(i)] W.h.p. it is possible to choose $n$ points
out of $n$ pairs of random points such that all components are smaller
than $f(c) n$.
\item[(ii)] W.h.p., for every choice of $n$ points out
of $n$ pairs of random points, the largest component has order at least
$g(c) n$.
\end{longlist}
Moreover, $f(c) \to0$ as $c\downarrow0$ and
$g(c) \to1$ as $c \to\infty$.
\end{theorem}

Again the result extends to the setting with $d\geq2$ choices; the
expression for $r$ then is $r=\sqrt[d+1]{c/n}$.

\subsection{Background and related work}
\label{sec:Background}

The notion of the ``power of choices'' essentially dates back to a 1994
paper by Azar, Broder, Karlin and Upfal \cite{AzarEtAl-STOC,AzarEtAl}.
In informal computer science terms, their result states that if one
allocates a large number of jobs to a large number of servers by
assigning each job to the currently less busy of \emph{two} randomly
chosen servers, one observes a dramatic improvement in load balancing
over a completely random assignment. This result marked the beginning
of the
development of the power of choices as a powerful new paradigm in
computer science, with
applications to load balancing, hashing, distributed computing, network
routing and other areas
(see \cite{SurveyPOC} for a comprehensive survey).

The mathematical model for this setting is usually given in terms of
balls and bins. In the standard balls and bins experiment, there are
$n$ balls and $n$ bins, and each
ball is dropped into a random bin (chosen uniformly at random,
independently of the choices for
the other balls). Denoting by $M_n$ the number of balls in the fullest
bin, also called the \emph{maximum load}, it is well known (and can be
proved by the first and second moment methods) that w.h.p. $M_n =
(1+o(1)) \ln n / \ln\ln n$ (an even more precise result is given
in \cite{Gonnet81}).

In the power of choices version of this setup, the $n$ balls arrive
sequentially, and for each ball \emph{two} random bins are sampled
(uniformly at random and independently from each other).
The goal now is to devise a strategy for choosing between the bins that
keeps the maximum load as small as possible.
An obvious choice of a strategy is the \emph{greedy} strategy where we
always choose the least full bin
(in case of a tie we can choose arbitrarily).
Azar et al. \cite{AzarEtAl} showed the following remarkable result.
Recall that a random variable $X$ \emph{stochastically dominates} the
random variable $Y$ if
$\Pee( X \geq x ) \geq\Pee( Y \geq x )$ for all $x\in\eR$.

\begin{theorem}[(\cite{AzarEtAl})]\label{thm:AzarEtAl}
Consider the power of choices balls and bins process with $n$ bins and
$n$ rounds, and
let $M_n$ denote the maximum load after the process ends if we employ
the greedy strategy.
Then
\[
M_n = \log\log n + O(1)\qquad \mbox{w.h.p. }
\]
Moreover, the maximum load under any other strategy stochastically dominates
the maximum load under the greedy strategy.
\end{theorem}

This result shows that in the power of two choices ball and bins process,
the maximum load under the greedy strategy is \emph{exponentially
smaller} than in the ordinary nonpower-of-choices version, and that
moreover it is very strongly concentrated.
It also shows that the greedy strategy is optimal in a very strong sense.

We should remark that Theorem~\ref{thm:AzarEtAl} is in a fact a slight
simplification of Theorem~1 of \cite{AzarEtAl}.
Among other things, it was also shown in \cite{AzarEtAl} that allowing
more than two choices per step further decreases the maximum load, but
only by a constant factor.
Note that the behaviour of the geometric power of choices process is in
contrast with this: As mentioned in the previous section, in our setting
a choice of $d$ options in each round results in a threshold at $r =
\Theta( \sqrt[d+1]{1/(n \cdot(\log\log n)^{d-1})} )$. Thus,
every additional choice per step increases the threshold by a power of $n$.

Theorem~\ref{thm:AzarEtAl} plays an important role in our proof of
Theorem~\ref{thm:online}, and may give some intuition why a power of
$\log\log n$ appears in the threshold for our geometric power of
choices process.

\textit{The Achlioptas process.}
The power of choices version of the classical Erd\H{o}s--R\'enyi graph
process is usually called the \emph{Achlioptas process} after Dimitris
Achlioptas, who first suggested it. The Achlioptas process starts with
an empty graph on $n$ vertices. In each round, two random vertex pairs
are presented, and the player needs to select exactly one of them for
inclusion as an edge in the evolving graph.
His goal is to delay or accelerate the occurrence of some monotone
graph property, such as containing a giant component, containing a
triangle, or containing a Hamilton cycle.

Bohman and Frieze \cite{Avoiding} were the first to study the
Achlioptas process. They showed that by an appropriate edge-selection
strategy, the emergence of a giant component can be delayed by a
constant factor. Several authors subsequently improved on their bounds,
and also showed that no improvement beyond a constant factor is
possible \cite{AvoidingUB,BirthControl}. The opposite problem of
creating a giant component as quickly as possible was studied in \cite
{Creating1,Creating2}, and an exact threshold for the offline problem
was determined in \cite{AchlioptasOfflineThreshold}. All these results
show that the power of choices offered by the Achlioptas process
affects the threshold for the appearance of a giant component only by a
constant factor.

More recently, the precise nature of the phase transition in the
Achlioptas process received much attention: countering ``strong
numerical evidence'' presented in \cite{AchlioptasScience}, Riordan and
Warnke \cite{AchlioptasContinuous} showed that for a large number of
natural player strategies, the Achlioptas phase transition is in fact
continuous.

Other properties that have been studied for the Achlioptas process
include Hamiltonicity \cite{AchlHamiltonicity} and the appearance of
copies of a given fixed graph $F$ \mbox{\cite{AchlSubgraphs1,AchlSubgraphs2,AchlCreatingSubgraphs}}.

\textit{The vertex Achlioptas process.}
The reader might wonder whether it is reasonable to compare our
geometric power of choices process to the Achlioptas process---after
all, we are selecting vertices in the former and edges in the latter.
Let us therefore consider the following process: the $n$ vertices of an
Erd\H{o}s--R\'{e}nyi random graph $G(n,m)$ (i.e., the random graph
sampled uniformly from all graphs on $n$ vertices and $m$ edges) are
revealed two at a time, along with all edges induced by the vertices
revealed so far, and we need to select one of the two vertices for
inclusion in a subgraph. Our goal is to avoid a linear-sized component
in the subgraph induced by the vertices we select.
This process indeed seems to be a better reference for our comparison,
but its phase transition has not been studied explicitly in the
literature. As it turns out, also in this {\em{``vertex Achlioptas
process''}}, w.h.p. a giant component cannot be avoided as soon as the
average degree of the underlying random graph exceeds a certain
constant. We give a proof for this in Appendix~\ref{sec:vertex-achlioptas}.

We remark that further discussion of our results and possible
directions for further work can be found
in Section~\ref{sec:concl}, at the end of the paper.

\subsection{About the proofs}

In the following, we informally outline some of the main ideas used in
the proofs of
Theorems \ref{thm:online} and \ref{thm:offline}. Our goal here is not
to give detailed proof sketches, but to give some impression of our
overall proof strategies and of the type of arguments used. We will
describe our proof strategies in more detail later where appropriate.

In all our proofs, we consider a discretized version of the random
geometric graph as follows: We divide the unit square into $\Theta
(r^{-2})$ many small squares (called \emph{boxes}) in such a way that,
essentially, we no longer need to worry about the precise locations of
the random points, but only need to know which of the boxes are
occupied by at least one point. In this way, our analysis of the
process reduces to an analysis of appropriate random subgraphs of a
large finite square grid (or king's move grid, see Section~\ref{sec:online.lb.pf} below). With some technical work, the results for
this grid then translate back to the original geometric setting to give
the desired results.

\begin{figure}

\includegraphics{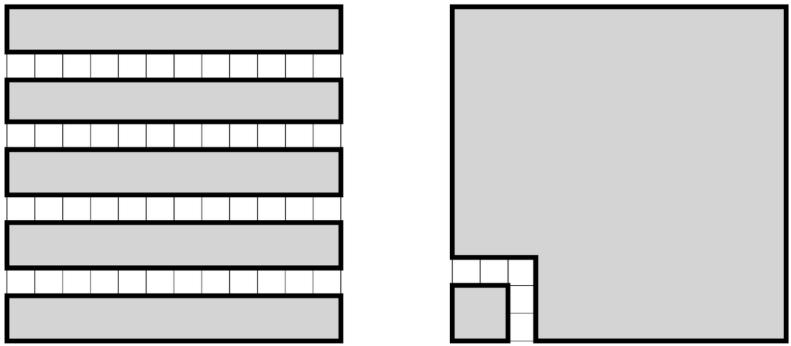}

\caption{Dividing the unit square using a barrier consisting of $K
\cdot(1/hr)$ blocks of dimensions $hr \times hr$,
when $K$ is large (left) and when $K$ is small (right).}\label{fig:moat}
\end{figure}

\textit{Lower bound proofs}. The key idea in the proofs of the lower
bound parts of both Theorems \ref{thm:online} and \ref{thm:offline} is
the following: To avoid the formation of large connected components of
occupied boxes (and thus also of large connected components in the
original power of choices random geometric graph setting), we designate
a subset of the boxes as the ``barrier''. This barrier separates the
unit square into ``small'' parts (see Figure~\ref{fig:moat}). The
player's goal is to prevent the appearance of paths of occupied boxes
crossing the barrier. In both settings, he selects points outside the
barrier whenever possible. His main worry are the pairs of points that
both fall into the barrier, as these force him to select a point in the barrier.

In the offline scenario, it is not too hard to show that even if the
barrier is only a (moderately large) constant number of boxes wide,
then w.h.p. there is a choice of points for which the barrier is not
crossed. The proof of this relies on an (approximate) analogy of our
setting with the Erd\H{o}s--R\'{e}nyi random graph in its subcritical
stage, where the boxes in the barrier play the role of the vertices and
the pairs of points that both fall into the barrier yield the edges.

In the online scenario, we can only ensure that the barrier will not be
crossed if it is at least $\Theta(\log\log n)$ boxes wide. Our strategy
here is more involved; the key part is to set things up in such a way
that, essentially, we can invoke Theorem~\ref{thm:AzarEtAl} to argue
that w.h.p. the barrier will not be crossed. To that end, we divide
the barrier into ``blocks'' consisting of $\Theta(\log\log n)\times
\Theta(\log\log n)$ boxes, and pay special attention to the blocks that
become ``dangerous'' because, oversimplifying slightly, only relatively
few more occupied boxes are required to create a (potentially crossing)
long path inside them.

\textit{Upper bound proofs}.
A key ingredient in both our upper bounds proofs is a simple
isoperimetric inequality for subgraphs of the square grid. It allows us
to conclude from the fact that there are relatively few unoccupied
boxes that a significant proportion of the occupied boxes must belong
to relatively large connected components of the graph induced by all
occupied boxes.

For the offline case, we use this as part of a combinatorial counting
argument. Essentially, we count the number of sets of boxes whose
removal decomposes the grid graph into ``small'' components. More
precisely, in order to keep the total number of sets to be considered
sufficiently low, we will only count appropriate \emph{subsets} of
decomposing sets as described. By an expectation argument, we will then
show that w.h.p. the player will not be able to avoid a single one of
these sets and, therefore, will be forced to create a ``large''
(linear-sized) component.

For the online case, we use a two-round approach and analyze the
process after $n/2$ points rounds as an intermediate step. To this end,
we divide the grid graph into ``b-blocks'' consisting of $O(1/(r\log
\log n))\times O(1/(r\log\log n))$ boxes, and define an appropriate
notion of ``good'' b-blocks [see Figure~\ref{fig:superblocks}(b) on page
\pageref{fig:good}]. We show that at time $n/2$, w.h.p. most boxes of
the original grid graph are occupied, and that as a consequence most
b-blocks are good. It follows with the mentioned isoperimetric
inequality that a bounded number of connected components of the graph
of b-blocks (defined in the obvious way---b-blocks are adjacent if
they share a side) covers most of the unit square.
Conditional on that, we then show that in the remaining $n/2$ rounds of
the process, w.h.p. at least one of these components will evolve into
a linear-sized component of the original grid graph. To show that the
player cannot avoid this, we apply a slight variation of Theorem~\ref
{thm:AzarEtAl}
in an appropriate geometric setup.


\section{Preliminaries}

For the sake of readability and clarity of exposition, we will mostly
ignore rounding, that is, we will usually omit floor and
ceiling signs. In all cases, it is a routine matter to check that all
computations and proofs also work
if floors and ceilings are added, and we leave this to the reader.

Throughout the paper, $\ln x$ will denote the \emph{natural logarithm} and
$\log x$ will denote the logarithm base $2$ (i.e., $\ln x$ is the
inverse of $e^x$ and $\log x$ is the inverse of~$2^x$).
For $n \in\eN$, we will denote $[n] := \{1,\ldots, n\}$.

By a slight abuse of notation, we will write $\Zed^2$ to denote the
\emph{graph of the integer lattice}, tht is, the
infinite graph with vertex set $\Zed^2$ and an edge between two
vertices if and only if their distance is exactly one.
Similarly, we will also identify $A \subseteq\Zed^2$ with the subgraph
of $\Zed^2$ it induces.
Thus, $[s]^2$ in particular denotes an $s\times s$ grid.

We will use the notation $\Bi(n,p)$ to denote the binomial distribution
with parameters $n$ and $p$; we
will use $\Po(\mu)$ to denote the Poisson distribution with mean~$\mu$;
and we will use
$\geom(p)$ to denote the geometric distribution with parameter~$p$.

We will use the following incarnation of the Chernoff bound.
A proof can, for instance, be found in \cite{PenroseBoek}, cf. Lemmas~1.1
and~1.2.

\begin{lemma}\label{lem:chernoff}
Let $X$ be a random variable with a binomial or Poisson distribution.
Then, letting $\mu:= \Ee X$, we have:
\begin{longlist}[(ii)]
\item[(i)] for $k \geq\mu$ we have
$\Pee( X \geq k ) \leq\exp[ - \mu H( k/\mu) ]$, and
\item[(ii)] for $k \leq\mu$ we have
$\Pee( X \leq k ) \leq\exp[ - \mu H( k/\mu) ]$,
\end{longlist}
where $H(x) := x\ln x - x + 1$.
\end{lemma}


\section{Lower bound proofs}
\label{sec:lb.pf}


\subsection{Proof of part \textup{(i)} of
Theorem \texorpdfstring{\protect\ref{thm:online}}{1}}
\label{sec:online.lb.pf}

We will consider a suitable discretization of the geometric graph.
Let $\Dscr_r$ be the obvious dissection into squares of side length~$r$,
that is,
%
\begin{equation}
\label{eq:Drdef} \Dscr_r := \bigl\{ \bigl[ir, (i+1)r\bigr]\times
\bigl[jr, (j+1)r\bigr] \dvtx 0 \leq i,j < 1/r \bigr\}.
\end{equation}
[We assume that $1/r$ is an integer throughout the section.\footnote{
The concerned reader may check that, if we take $\tilde{h} := \lfloor
100\log\log n\rfloor$ instead of $h
= 100\log\log n$ and
$\rtil:= 1 / \tilde{h}\lfloor1/r\tilde{h}\rfloor$ instead of $r$,
then $1/r$ and $1/hr$ are both integers, and
$\tilde{h} = (1+o(1)) h, \tilde{r} = (1+o(1)) r$, and all proofs and
computations
in this section carry through. Since $\rtil\geq r$, this also
establishes the result for $r$.
}]
We refer to the elements of $\Dscr_r$ as \emph{boxes}, and we
say that a box is \emph{occupied} if it contains a point of our
process and
\emph{empty} otherwise.

The \emph{king's move grid} $\King_s$
on $s\times s$ vertices is the graph with vertex set $[s]^2$ and an
edge between two vertices if and only if their distance is at most
$\sqrt{2}$.
This way, we can also move diagonally, which explains the name ``king's
move grid''---at least to those familiar with the rules of chess.

We will identify the boxes of the dissection $\Dscr_r$ with the
vertices of the king's move grid $\King_{(1/r)}$---that is, we
consider two boxes adjacent if they share a side or a corner. Note that
if two points are adjacent in the original geometric graph, then they
must lie in boxes that are adjacent in $\King_{(1/r)}$. Thus, denoting
by $\Occ_r$ the subgraph of $\King_{(1/r)}$ induced by the occupied
boxes, the following holds: If $C_1, C_2$ are distinct components of
$\Occ_r$, then the points in $C_1$ and the points in $C_2$ belong to
different components of the geometric graph.

We further group the boxes into ``blocks'' consisting of $h \times h$
boxes each, where $h = 100 \cdot\log\log n$.
(For convenience, we assume that $1/r$ and $1/hr$ are both integers.)
Again it is useful to consider the blocks as vertices of a king's move
grid $\King_{(1/hr)}$.

Our strategy can be described as follows.
At the start of the game, we will pick a constant $K = K(c)$, to be
made explicit later on, and we pick a set of
\[
N := K \cdot(1/hr),
\]
blocks that forms a ``barrier'' as in Figure~\ref{fig:moat} that we will
``defend''. The next lemma formally captures the essential properties of
our ``barrier''. Its proof is indicated in Figure~\ref{fig:moat}; the
details are left to the reader.

\begin{lemma}\label{lem:moat}
For every fixed $K > 0$, the following holds for all sufficiently large
$s \in\eN$.
There is a subgraph $H \subseteq\King_s$ of the king's move grid with
$v(H) \geq s^2 - Ks$ vertices such that
\[
v( H_{\max} ) \leq\bigl(a(K)+o_s(1)\bigr) \cdot
s^2,
\]
where $H_{\max}$ denotes the largest component of $H$, and
\[
a(K) = \cases{1 / \bigl(\lfloor K\rfloor+ 1 \bigr), &\quad
$\mbox{if } K > 1,$
\vspace*{2pt}\cr
1-(K/2)^2, &\quad $\mbox{if } K \leq1.$}
\]
In particular, $0 < a(K) < 1$ for all $K > 0$ and
$a(K) \to0$ as $K \to\infty$. 
\end{lemma}

Let us remark that while it is clearly possible to improve on the expression
for $a(K)$ given in Lemma~\ref{lem:moat} we have made no attempt to do so
as the current version of the lemma suffices for our purposes.

Hence, our strategy begins by selecting $N$ blocks according to
Lemma~\ref{lem:moat}
so that they form a barrier $\Moat$ that divides the unit square
$[0,1]^2$ into pieces of
area at most $a(K)$.
During the game, we will now attempt to prevent the formation of
(king's move) paths of occupied boxes crossing the barrier.
We will devise a strategy that will succeed in defending the barrier in
this sense (w.h.p.).
This then implies that the area of the occupied boxes corresponding to
the largest component
of the geometric graph is $a(K)+o(1)$. [This may include some boxes
from the barrier; note, however, that the entire barrier has area $o(1)$.]
With an additional argument showing that (w.h.p.) no union of boxes
with area $A\geq\eps$ ($\eps>0$ arbitrary but fixed) will have
received more than $(1+\eps)\cdot A\cdot n$ points, it then follows
that the number of vertices of the largest component of the geometric
graph is bounded by $(a(K)+o(1)) \cdot n$, which completes the proof.

\begin{figure}

\includegraphics{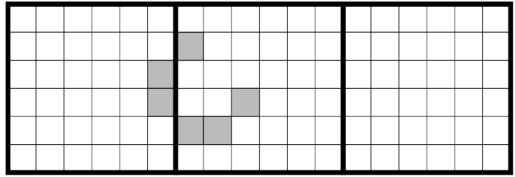}

\caption{The left and middle blocks are bad. (The gray cells represent
occupied cells.)}\label{fig:badblock}
\end{figure}

Let us now describe our strategy in more detail.
We say that a block $B$ is \emph{bad} if there is a path of occupied
boxes of length $h$ that uses
one of the boxes of $B$ and uses only boxes of the barrier.
Observe that such a path may use boxes in blocks other than $B$, but it
may not use
boxes in blocks that are either not adjacent to $B$ or do not belong to
the barrier.
See Figure~\ref{fig:badblock} for a depiction.

Clearly, if there are no bad blocks then there cannot be any path
crossing the barrier.

We will say that a block $B$ is \emph{dangerous} if occupying up to $2
\log\log n$ additional boxes
can render it bad.
In other words, $B$ is dangerous if there is some path $P$ of length
$h$ that uses only boxes of the barrier, and
at least one box of $B$, such that at least $h - 2\log\log n$ boxes of
$P$ are occupied.
If $B$ is not dangerous, we will call it \emph{safe}.
Let $N_{\mathrm{bad}}(t), N_{\mathrm{dang}}(t)$ denote the number of
bad respectively dangerous blocks at the end of round $t=1,\ldots,n$.

We will keep an ordered list of blocks $L(t) = (L_1(t), \ldots,
L_{2^{-h}N}(t) )$
containing exactly $2^{-h}N$ blocks, which will be updated at the end
of each round.
We will make sure that if $\Ndang(t) \leq2^{-h} N$ then $L(t)$
contains all dangerous blocks by applying
the following update rule.
If $\Ndang(t) < 2^{-h} N$ and some previously safe block $B \notin
L(t)$ becomes dangerous
during round $t+1$, then we replace an arbitrary safe block in the list
with the new dangerous block.
[If for some $t$ it happens that $\Ndang(t) > 2^{-h} N - k$ and $k$
blocks become dangerous during
round $t+1$, then our strategy will have failed.]

We will call the blocks in $L(t)$ \emph{pseudodangerous}, and the blocks
not in $L(t)$
\emph{pseudosafe} (with respect to round $t$).

Our strategy can now be described as follows:
\begin{longlist}[\strat{1}]
\item[\strat{1}] We always pick a point outside of the barrier $M$ if
we can. If both points are outside the barrier, we choose randomly.
\item[\strat{2}] If both points fall inside the barrier and both are in
pseudosafe blocks, then we play randomly.
\item[\strat{3}] If both fall inside the barrier, one in a
pseudodangerous, and one in a pseudosafe block, then we choose
the pseudosafe block.
\item[\strat{4}] If both fall in the barrier, both in pseudodangeous
blocks, say in $L_i(t)$ and $L_j(t)$,
then we do the following.
Set
\[
k_i := \bigl| \bigl\{ t' < t \dvtx \mbox{ we played in
$L_i\bigl(t'\bigr)$ in round $t'$} \bigr
\}\bigr |,
\]
and define $k_j$ similarly.
If $k_i < k_j$, then we play in $L_i(t)$ and if $k_i > k_j$ then we
play in $L_j(t)$.
In case of a draw, we play randomly.
\end{longlist}

\begin{remark*} Recall that which block $L_i(t)$ points to may change
between rounds.
Let us stress that in \strat{4} we compare the number of times we
played in the $i$th, respectively, the $j$th pseudorandom block as opposed
to the number of times we played in the specific blocks that $L_i(t)$
respectively  $L_j(t)$ represent in round $t$.
On the other hand, once $L_i(t)$ points to a dangerous block, the value
of $L_i(t')$ will remain the
same for all $t' \geq t$. Thus, if $L_i(t)$ points to a dangerous block~$B$,
the number of times we played in the $i$th pseudorandom block until
round $t$ is an upper bound on
the number of times we played in $B$ \emph{since it became dangerous}.
This subtle, but very important point shows that---provided the number
of dangerous blocks remains
below $2^{-h}N$---if the number of rounds in which we play in the
$i$th pseudorandom block stays below $2\log\log n$ for
every index $i$ then no bad blocks will be formed.
\end{remark*}

As indicated before, we will show that our strategy works by eventually
proving that, w.h.p., no bad blocks will be formed (which
implies the barrier will not get crossed), and then showing that
for each of the regions that the barrier divides the unit square into,
the number of points that fell into the
region is proportional to its area (w.h.p).

We start by proving the following lemma which shows that, w.h.p., we
never have more
dangerous blocks than entries in the list $L(t)$. Clearly, this implies
that we will be able to keep the dangerous blocks a subset of the
pseudodangerous blocks.

\begin{lemma}\label{lem:Dsmall}
There is an absolute constant $C_0 > 0$ such that if $K \leq C_0 / c$,
then $\Ndang(t) \leq2^{-h} N$ for all $t=1,\ldots, n$, w.h.p.
\end{lemma}

\begin{pf}
Observe that the area of the barrier is
\[
\area(\Moat) := N \cdot(hr)^2 = K h r.
\]

Let $R$ denote the number of rounds in which we (are forced to) play
inside the
barrier $\Moat$. Clearly, $R \isd\Bi( n, \area^2(\Moat) )$.
Note that
%
\begin{equation}
\label{eq:ER} \Ee R = n \cdot K^2 \cdot h^2 \cdot
r^2 = n^{{1}/3+o(1)}
\end{equation}
by choice of $r$.

To facilitate the analysis, it is helpful to consider what would happen
under a slightly different setup.
Suppose that whenever both points fall inside the barrier during some
round we add both points, and otherwise
we pick a point outside the barrier.
Clearly, we will end up with $2R$ points in the barrier, distributed uniformly.
Moreover, the set of occupied boxes inside the barrier under this setup
will be a superset of the set of
occupied boxes in the original setup. 

Let $Z \isd\Po( 8 \Ee R )$ be a Poisson variable with mean $8 \Ee R$.
By applying the Chernoff bound, Lemma~\ref{lem:chernoff}, together
with \eqref{eq:ER} we can easily see that
\begin{eqnarray*} \Pee( Z < 2R ) & \leq& \Pee( Z \leq\Ee Z / 2 ) +
\Pee( R \geq2\Ee R )
\\
& \leq& \exp\bigl[ - \Ee Z \cdot H\bigl(\tfrac{1}2\bigr) \bigr] +
\exp\bigl[ - \Ee R \cdot H(2) \bigr]
\\
& = & \exp\bigl[ - \Omega\bigl( n^{{1}/3+o(1)} \bigr) \bigr].
\end{eqnarray*}

Let us once more modify the setup slightly, and just drop $Z$ points on
the barrier (their locations
chosen uniformly at random and independent of $Z$ and the locations of
the other points).
This way, the points in the barrier will form a \emph{Poisson process},
which has the convenient consequence that
the events that different boxes in the barrier are occupied become
independent (see, e.g., \cite{KingmanBoek}).
Note that, by choosing a suitable coupling, we can ensure that if
$Z\geq2R$ then this setup dominates our previous setup in the sense that
the set of points in the barrier under the old setup is a subset of the
set of points in the barrier
under the new setup. This allows us to bound the probability that a
block becomes dangerous by considering the new setup.
Observe that the expected number of points in a given box is
\begin{eqnarray*} \mu & := & \Ee Z \cdot \biggl(\frac{r^2}{\area(\Moat)}
\biggr)
\\
& = & 8 \cdot n \cdot K^2 \cdot h^2 \cdot
r^2 \cdot \biggl(\frac{r^2}{K
\cdot
h \cdot r} \biggr)
\\
& = & 8 \cdot K \cdot h \cdot n \cdot r^3
\\
& = & 800 \cdot K \cdot c.
 \end{eqnarray*}
Let us write $p := 1 - e^{-\mu}$.
Then a box in the barrier is occupied with probability~$p$,
independently of all other boxes.

Next, observe that the number of paths of length $h$ starting either in
a given block or in one of
the neighbouring blocks is at most
$9 h^2 8^{h-1}$ (such a path starts in one of the $9 h^2$ boxes
belonging to the block and
adjacent blocks, and there are always at most 8 choices for the next
box of the path).
Let $\pdang$ denote the probability that a particular block is
dangerous under this new setup.
The union bound gives
\[
\pdang \leq 9 \cdot h^2 \cdot8^{h-1} \cdot\Pee \bigl( \Bi(
h, p ) \geq h - 2\log \log n \bigr).
\]
Using the Chernoff bound (Lemma~\ref{lem:chernoff}), we see that for
all $k\geq h$
\begin{eqnarray*}\Pee \bigl( \Bi( h, p ) \geq h - 2\log\log n
\bigr) & \leq& \exp \biggl[ - hp \cdot H \biggl(\frac{h-2\log\log n}{hp} \biggr) \biggr]
\\
& = & \exp \biggl[ - hp \cdot H\biggl(\frac{98}{100 \cdot p}\biggr) \biggr],
\end{eqnarray*}
where as usual $H(x) = x\ln x - x + 1$.

Observe that $p H ( \frac{98}{100p} ) \to\infty$ as
$p\downarrow0$, since
\begin{eqnarray*}
p H \biggl( \frac{98}{100 \cdot p} \biggr) &=& p \cdot \biggl( \biggl(\frac{98}{100p}
\biggr)\ln \biggl(\frac
{98}{100p} \biggr)-\frac{98}{100p}+1 \biggr) \\
&=&
\frac{98}{100}\ln \biggl(\frac{98}{100p} \biggr) + \frac{98}{100} + p.
\end{eqnarray*}
Hence, there exists a universal constant $p_0$ such that
%
\begin{equation}
\label{eq:pineq} \exp \biggl[ - p H \biggl(\frac{98}{100 \cdot p} \biggr) \biggr] \leq
\frac
{1}{100}\qquad \mbox{for all } p\leq p_0.
\end{equation}
Let us set $C_0 := - \ln(1-p_0) / 800$, so that
$K = K(c) \leq C_0 / c$ implies that
$p = 1 - \exp[ - 800 K c ] \leq p_0$.

With this choice of $K$, we have
\begin{eqnarray*}\pdang & \leq& 9 \cdot h^2
\cdot8^{h-1} \cdot\exp \biggl[ - hp \cdot H\biggl(\frac{98}{100
\cdot p}
\biggr) \biggr]
\\
& = & \frac{9}{8} \cdot h^2 \cdot \biggl( 8 \cdot\exp\biggl[
- p H\biggl(\frac{98}{100
\cdot p}\biggr) \biggr] \biggr)^h
\\
& \leq& \frac{9}{8} \cdot h^2 \biggl( \frac{8}{100}
\biggr)^h
\\
& \leq& 10^{-h},
\end{eqnarray*}
where the last line holds for $n$ sufficiently large (recall $h\to
\infty
$ as $n\to\infty$).

Let $D$ denote the number of dangerous blocks in our modified process
where we simply
drop $Z$ points uniformly at random on the barrier.
Then $\Ee D = N \cdot\pdang$, and hence
\begin{eqnarray*} \Pee\bigl( D > 2^{-h} N \bigr) & = &
\Pee\bigl( D > \bigl(2^{-h}/\pdang\bigr) \cdot\Ee D \bigr)
\\
& \leq& \Pee\bigl( D > 5^h \cdot\Ee D \bigr)
\\
& \leq& 5^{-h}
\\
& = & o(1),
\end{eqnarray*}
where we used Markov's inequality for the third line.

Putting everything together, we see that
\begin{eqnarray*} \Pee\bigl( \Ndang> 2^{-h} N \bigr) &
\leq& \Pee( Z < 2R ) + \Pee\bigl( D > 2^{-h} N \bigr)
\\
& = & o(1),
\end{eqnarray*}
which completes the proof of the lemma.
\end{pf}

Let $\Rpsd$ denote the number of rounds in which we (are forced to) play
in a pseudodangerous block.

\begin{lemma}\label{lem:Msmall}
$\Rpsd\leq2^{-h} N$ w.h.p.
\end{lemma}

\begin{pf}
The probability that both points land in pseudodangerous blocks
in round $t=1,\ldots, n$ is
\[
p := \bigl( 2^{-h} N h^2 r^2
\bigr)^2. 
\]
Thus, we have
\begin{eqnarray*} \Ee\Rpsd/ \bigl(2^{-h} N\bigr) & = &
(np) / \bigl(2^{-h} N\bigr)
\\
& = & \bigl(n 2^{-2h} N^2 h^4 r^4
\bigr) / \bigl(2^{-h} N\bigr)
\\
& = & n 2^{-h} N h^4 r^4
\\
& = & O \bigl( n 2^{-h} h^3 r^3 \bigr)
\\
& = & O \bigl( h^2 2^{-h} \bigr)
\\
& = & o(1),
\end{eqnarray*}
where we used that $N = O( 1/hr )$ in the fourth line, and that $nr^3 =
c / \break \log\log n = \Theta( 1/h )$
in the fifth line.
It follows that
\[
\Pee\bigl( \Rpsd> 2^{-h} N \bigr) \leq \frac{\Ee\Rpsd}{2^{-h} N} = o(1),
\]
by Markov's inequality.
\end{pf}

\begin{lemma}\label{lem:Nbadisnul}
If $K \leq C_0 / c$ with $C_0$ as in Lemma~\ref{lem:Dsmall}, then
$\Nbad
(n) = 0$ w.h.p.
\end{lemma}

\begin{pf}
If some bad block got created even though we stuck to our strategy,
then it must be the case that either (a) there were more
dangerous blocks than places in our list of pseudodangerous blocks, or
(b) there is some index $1\leq i \leq2^{-h}N$ such that there were
more than $2 \log\log n$ rounds $t$ when we played the $i$th
pseudodangerous block in our list.
[Recall the remark immediately following the description of the
strategy \strat{1}--\strat{4}.]

Let $E$ denote the event that (b) happens.
Clearly,
\begin{eqnarray*}
&&\Pee\bigl( \Nbad(n) > 0 \bigr) \\
&&\qquad\leq \Pee\bigl( \Ndang(n) > 2^{-h} N
\bigr) + \Pee\bigl( \Rpsd> 2^{-h} N \bigr) + \Pee\bigl( E \mbox { and }
\Rpsd\leq2^{-h} N \bigr).
\end{eqnarray*}
By Lemmas \ref{lem:Dsmall} and \ref{lem:Msmall}, the first two terms of
the right-hand side
are $o(1)$.
Now notice that
\begin{eqnarray*}\Pee\bigl( E \mbox{ and } \Rpsd
\leq2^{-h} N \bigr) & \leq& \Pee\bigl( E | \Rpsd\leq2^{h} N
\bigr)
\\
& \leq& \Pee\bigl( E | \Rpsd= 2^{-h} N \bigr),
\end{eqnarray*}
where the second inequality holds by obvious monotonicity.
Now notice that, by part \strat{4} of our strategy,
the event that $E$ holds, given that $\Rpsd= 2^{-h} N$, can be viewed
as the event that, in the standard power of choices balls and bins
setup with
$\tilde{n} := 2^{-h} N$ balls and $\tilde{n}$ bins, the maximum load is
at least $2 \log\log n$.
Since $\tilde{n}= 2^{-h} \cdot K \cdot(1/hr) = n^{{1}/3 + o(1)}$ we
have that
\[
\log\log\tilde{n} = \bigl(1+o(1)\bigr) \log\log n.
\]
It therefore follows immediately from Theorem~\ref{thm:AzarEtAl}
that $\Pee( E | \Rpsd= 2^{-h} N ) = o(1)$.
We see that $\Pee( \Nbad(n) > 0 ) = o(1)$, as required.
\end{pf}

By this last lemma, our strategy succeeds (w.h.p.) in confining the
components of the
evolving random geometric graph to subsets of the unit square of area
bounded by
$a(K)+o(1)$.
The finishing touch of the proof of part (i) of
Theorem~\ref{thm:online} comes in the form of the following lemma.

\begin{lemma}\label{lem:online.lb.area}
For every $\eps> 0$, the following holds if we follow the strategy set
out above.
W.h.p. every $A \subseteq[0,1]^2$ that is the union of boxes of the
dissection $\Dscr_r$
and with $\area(A) \geq\eps$ contains at most $(1+\eps)\cdot\area(A)
\cdot n$ points.
\end{lemma}

\begin{pf}
For $A \subseteq[0,1]^2$, let $\Nscr(A)$ denote the number of points
in $A$
(in round~$n$).
Let us first recall that the barrier satisfies $\area(\Moat) = K h r
= o(1)$.
If $R$ denotes the number of rounds in which we (are forced to) take a
point in the barrier then clearly
$R \isd\Bi(n, \area^2(\Moat) )$, and in particular
$\Ee R = n \cdot\area^2(\Moat) = o(n)$.
Hence, by Markov's inequality
\[
\Pee\bigl( \Nscr(\Moat) \geq\bigl(\eps^2/2\bigr) n \bigr) = \Pee
\bigl( R \geq\bigl(\eps^2/2\bigr) n \bigr) = o(1).
\]
Let $\Ascr$ denote all the subsets of $[0,1]^2$ under consideration,
that is, all
$A$ that are unions of boxes of $\Dscr_r$ and have area at least $\eps$.

Pick an arbitrary $A \in\Ascr$ and set $A' := A \setminus\Mscr$.
Then $\area(A') = \area(A) - o(1)$. Note that in every round $1 \leq t
\leq n$, a
point is added to $A'$ with probability
%
\[
p := \bigl(1-\area^2(\Moat)\bigr) \cdot\frac{\area(A')}{1-\area(\Moat)} =
\bigl(1+o(1)\bigr)\cdot\area(A).
\]
(This is because if both points fall in the barrier, we obviously add a
point outside of
$A'$, and otherwise we add a point drawn according to the uniform
distribution on
$[0,1]^2 \setminus\Moat$.)
We have $\Nscr(A') \isd\Bi(n,p)$ so that by the Chernoff bound
(Lemma~\ref{lem:chernoff})
we have
\begin{eqnarray*} \Pee\bigl( \Nscr\bigl(A'\bigr) >
(1+\eps/2) \cdot\area(A) \cdot n \bigr) & \leq& \exp\biggl[ - n \cdot p \cdot H
\biggl(\frac{(1+\eps/2) \cdot\area(A)
\cdot
n}{n p} \biggr) \biggr]
\\
& = & \exp\bigl[ - \Omega( n ) \bigr],
\end{eqnarray*}
since $p = \Omega(1)$ and $(1+\eps/2) \cdot\area(A) \cdot n / (n p) =
1 + \eps/2 + o(1)$ is bounded away from one.
Let $E$ denote the event that there exists $A\in\Ascr$ with $\Nscr
(A) >
(1+\eps) \area(A) n$. Noting that $(\eps^2/2) n\leq(\eps/2)\cdot
\area
(A)\cdot n$ for all $A\in\Ascr$, we obtain that
\begin{eqnarray*}\Pee(E) & \leq& \Pee\bigl( \Nscr(\Moat) > \bigl(
\eps^2/2\bigr) n \bigr) + \sum_{A \in\Ascr} \Pee
\bigl( \Nscr(A \setminus\Moat) > (1+\eps/2) \cdot\area (A) \cdot n \bigr)
\\
& \leq& o(1) + 2^{(1/r)^2} \cdot\exp\bigl[ - \Omega( n ) \bigr]
\\
& = & o(1) + \exp\bigl[n^{{2}/3+o(1)} - \Omega(n) \bigr]
\\
& = & o(1),
\end{eqnarray*}
where in the second line we used that $|\Ascr| \leq2^{(1/r)^2}$ as
sets in $\Ascr$ are unions of the boxes of our dissection $\Dscr_r$,
and in the third line we used the specific form of $r$.
\end{pf}

Lemmas \ref{lem:Dsmall}--\ref{lem:online.lb.area} together imply
part (i) of Theorem~\ref{thm:online}.
For completeness, we spell out the details.

\textit{Proof of part \textup{(i)} of
Theorem~\ref{thm:online}}:
If we take $K = C_0 / c$ with $C_0$ as provided by Lemma~\ref
{lem:Dsmall}, and follow the strategy described by \strat{1}--\strat
{4} above,
then by Lemma~\ref{lem:Nbadisnul}, w.h.p. every connected component of
the resulting geometric graph
will lie inside a set of boxes of area at most $a(K)+o(1)$ with $a(\cdot)$
as in Lemma~\ref{lem:moat}.
Therefore, by Lemma~\ref{lem:online.lb.area}, w.h.p., every component
of the geometric graph
will have at most $n \cdot(a(K)+o(1))$ vertices. Thus, the claim
follows for, say, $f(c) := \sqrt{a(C_0/c)}$. [Recall that $0<a(K)<1$
for all $K$.] Since $K=C_0/c \to\infty$ as $c\downarrow0$ and
$a(K)\to0$ as $K\to\infty$, we also have $f(c) \to0$ as $c
\downarrow0$.


\subsection{Proof of part \textup{(i)} of
Theorem \texorpdfstring{\protect\ref{thm:offline}}{2}}
\label{sec:offline.lb.pf}

Our proof strategy is similar to the one for part (i)
of Theorem~\ref{thm:online} used in the preceding section. We will make
use of a standard result for the Erd\H{o}s--R\'enyi random graph
$G(n,m)$. Recall that $G(n,m)$ is obtained by taking a set of $n$
vertices, and selecting a set of $m$ edges uniformly at random from all
possible sets of $m$ edges. A graph is 1-\emph{orientable} if its edges
can be oriented in such a way that every vertex has indegree at most $1$.

The following result is a special case of Theorem~5.5 in the standard
reference~\cite{randomgraphs}.

\begin{theorem}\label{thm:erdosrenyi}
If $m \leq cn$ with $c < \frac{1}2$ then $G(n,m)$
consists only of trees and unicyclic components, w.h.p. In particular,
$G(n,m)$ is $1$-orientable w.h.p.
\end{theorem}

Let $\Dcal_r$ again be defined by \eqref{eq:Drdef}.
Again we will consider blocks, that is, $h\times h$ groups of boxes.
However, this time we simply set
$h := 100$.
Again we pick a constant $K=K(c)$, and build a barrier $\Mcal$
consisting of $N := K \cdot(1/hr)$ blocks,
in such a way that the barrier divides the unit square into parts of
area no more than $a(K)$ with $a(\cdot)$ as in Lemma~\ref{lem:moat}. As
before, we select points outside the barrier whenever possible,
breaking ties randomly. As we will see, Theorem~\ref{thm:erdosrenyi}
will then allow us to deal relatively easily with pairs of points that
both fall into the barrier.
Let $R$ denote the number of rounds in which both points fall into the barrier.

\begin{lemma}\label{lem:Roffline}
There is an absolute constant $C_0$ such that if $K \leq C_0 / c$ then
$R < N/100$ w.h.p.
\end{lemma}

\begin{pf}
Observe that $R \isd\Bi( n, \area^2(\Mcal) )$, and
$\area(\Mcal) = N \cdot(hr)^2 = Khr$.
Hence,
\[
\Ee R = n K^2 h^2 r^2 = \frac{c K^2 h^2}{r},
\]
using that $nr^3 = c$. This is less than $N/100 = K/(100hr)$ for
$K<1/(100h^3c)=1/(10^8c)$. Hence, the claim follows for $C_0:=10^{-8}$.

For such a choice of $K$, the Chernoff bound (Lemma~\ref
{lem:chernoff}) gives
\begin{eqnarray*} \Pee( R > N / 100 ) & \leq& \Pee( R > 10 \cdot\Ee
R )
\\
& \leq& \exp\bigl[ - \Omega( \Ee R ) \bigr]
\\
& = & \exp\bigl[ -\Omega\bigl( n^{1/3} \bigr) \bigr]
\\
& = & o(1),
\end{eqnarray*}
where we used the fact that
$\Ee R = \Theta( 1/r) = \Theta( n^{1/3} )$.
\end{pf}

Let us say that a block gets \emph{doubly hit} in some round if both
balls fall into the block in that round.

\begin{lemma}\label{lem:doublyhit}
W.h.p., no block of the barrier gets doubly hit in more than three rounds.
\end{lemma}

\begin{pf}
Let us fix a block $B$, and let $Z$ denote the number of rounds in
which $B$ gets doubly hit.
Clearly, $Z \isd\Bi( n, (hr)^4 )$.
Then we have
\[
\Pee( Z \geq4 ) = O\bigl( n^4 r^{16} \bigr) = O\bigl(
n^{-4/3} \bigr).
\]
Thus, the expected number of blocks that get doubly hit in at least
four rounds is
$O( N \cdot n^{-4/3} ) = O( n^{-1} ) = o(1)$.
\end{pf}

Let us now define an auxiliary (random) graph $\tilde{G}$, whose
vertices are the $N$
blocks of the barrier and where for every pair of blocks $B\neq B'$
there is an edge between
them if in some round one of the points landed in $B$ while the other
landed in $B'$.

\begin{lemma}\label{lem:Gtil}
Provided $K \leq C_0/c$ with $C_0$ as in Lemma~\ref{lem:Roffline},
the graph $\tilde{G}$ is 1-orientable, w.h.p.
\end{lemma}

\begin{pf}
Let $E$ denote the event that $\tilde{G}$ consists of trees and
unicyclic components (and hence is $1$-orientable), and
let $R'$ denote the number of rounds in which the points fell into two
different blocks of the barrier.
Observe that if we condition on $R'=m$ then $\tilde{G}$ is just a
copy of the Erd\H{o}s--R\'enyi random graph $G(N, m)$.
Let $F$ denote the event that the Erd\H{o}s--R\'enyi random graph
$G(N, N/100)$ consists of trees and unicyclic components.

Then we have that
\[
\Pee\bigl( E^c \bigr) \leq\Pee\bigl( R' > N/100 \bigr)
+ \Pee\bigl( F^c \bigr) = o(1),
\]
by Lemma~\ref{lem:Roffline} and Theorem~\ref{thm:erdosrenyi}.
\end{pf}

As mentioned, our strategy for the offline process will always select a
point outside the barrier if possible, choosing randomly if both points
fall outside the barrier. (Which point we select when both points
fall in the barrier will be specified shortly.) Let us note that the
proof of Lemma~\ref{lem:online.lb.area} in the previous section only
used part \strat{1} of our strategy for the online setting and,
therefore, carries over to our offline strategy.

\begin{lemma}\label{lem:offline.lb.area}
Consider the offline process, with the dissection $\Dcal_r$, the
barrier $\Mcal$, etc., as above. Assume that we
always select a point outside the barrier if we can, choosing randomly
if both points
fall outside the barrier.
Then for every $\eps> 0$ the following holds w.h.p.:
every $A \subseteq[0,1]^2$ that is the union of boxes of the
dissection $\Dscr_r$
and with $\area(A) \geq\eps$ contains at most $(1+\eps)\cdot\area(A)
\cdot n$ points. 
\end{lemma}

We are now ready for the proof of part (i) of
Theorem~\ref{thm:offline}.

\textit{Proof of part \textup{(i)} of
Theorem~\ref{thm:offline}}:
We show that, w.h.p., we can select one point from each pair in such a
way that no block of the barrier will contain more than four points.
Clearly, this then implies that the barrier will not be crossed.

To see this, note first that, by Lemma~\ref{lem:doublyhit}, w.h.p., no
block will contain more
than three points coming from rounds when it was doubly hit. Observe
also that, by Lemma~\ref{lem:Gtil}, w.h.p., the auxiliary graph
$\tilde
{G}$ is $1$-orientable. Hence, it is possible to select one point from
each pair of points that both fall into the barrier in such a way that
this contributes at most one point to each block. Thus, we can ensure
that, in total, each block will indeed contain at most four points (w.h.p.).

This shows that the player succeeds (w.h.p.) in stopping the barrier
from getting crossed.
The result now follows with Lemma~\ref{lem:offline.lb.area} exactly as
in the online case.


\section{Upper bound proofs}
\label{sec:ub.pf}

\subsection{Preliminaries}

We need a number of auxiliary results for our upper bound proofs. We
collect these in the next few subsections.
\subsubsection{Isoperimetric inequalities}

Recall that we identify subsets of $\mathbb{Z}^2$ or $[s]^2$ with the
subgraphs of the infinite integer grid induced by them.

\begin{lemma}\label{lem:isoperimetric}
Suppose $H \subseteq\Zed^2$ is a finite induced subgraph of the
integer lattice.
Then
\[
e\bigl( H, H^c \bigr) \geq4 \sqrt{v(H)}.
\]
\end{lemma}

\begin{pf}
Let $H$ be an arbitrary finite induced subgraph of $\mathbb{Z}^2$,
let $H_x$ denote the projection of $H$ on the $x$-axis, and let $H_y$
denote the
projection on the $y$-axis.
Let us write $\ell_x := |H_x|, \ell_y = |H_y|$.
Note that every vertical line that intersects $H$ contributes at least
two vertical edges to $e(H, H^c)$,
and that analogously every horizontal line that intersects $H$
contributes at least two horizontal edges. Thus,
\[
e\bigl(H, H^c \bigr) \geq2 \ell_x + 2
\ell_y.
\]
On the other hand, it is clear that
\[
v(H) \leq\ell_x \cdot\ell_y.
\]
Thus, we obtain with straightforward calculus that
%
\begin{equation}
\label{eq:vHmaxeq} v(H) \leq\max_{0 \leq z \leq{e(H,H^c)}/{2}} z \cdot \biggl(
\frac
{e(H,H^c)}{2}-z \biggr) = \biggl(\frac{e(H, H^c)}{4} \biggr)^2,
\end{equation}
which is equivalent to the claim.
\end{pf}

We will need the following strengthening of Lemma~\ref{lem:isoperimetric}.

\begin{lemma} \label{lem:new}
Let $x>0$ be given, and suppose $H \subseteq\Zed^2$ is a finite
induced subgraph of the
integer lattice such that all connected components of $H$ have at most
$x$ vertices.
Then
\[
e\bigl( H, H^c \bigr) \geq4 \frac{v(H)}{\sqrt{x}}.
\]
\end{lemma}

\begin{pf}
Let $\mathcal{C}$ denote the set of components of $H$. By Lemma~\ref
{lem:isoperimetric}, each component $C$ satisfies $\sqrt{v(C)}\leq
\frac
{e(C, C^c)}{4}$. Hence,
\[
v(H)=\sum_{C\in\mathcal{C}}v(C) \leq\sqrt{x}\sum
_{C\in\mathcal{C}}\sqrt{v(C)} \leq\sqrt{x}\sum
_{C\in\mathcal{C}} \frac{e(C, C^c)}{4} = \sqrt{x}\cdot\frac{e(H, H^c)}{4},
\]
which is equivalent to the claim.
\end{pf}

\begin{lemma}\label{lem:ExtremalMother}
Let $\alpha, \beta> 0$ and $s \in\eN$ be given.
Suppose $H \subseteq[s]^2$ is an induced subgraph of the $s\times s$
grid with $v( H ) \geq s^2 - \alpha s$. Moreover, let $H_\beta
\subseteq H$ denote
\[
H_\beta:= \cup\bigl\{ C \subseteq H \dvtx C \mbox{ is a component of
$H$, and } v(C) \leq\beta s^2 \bigr\},
\]
that is, $H_\beta$ is the union of all components of $H$ with at
most $\beta s^2$ vertices.
Then $v(H_\beta) \leq\sqrt{\beta} \cdot(1+\alpha) \cdot s^2$.
\end{lemma}

\begin{pf}
Let $H^c$ denote the complement of $H$ in $\mathbb{Z}^2$ (not $[s]^2$),
and observe that every edge of $e(H, H^c)$ connects a vertex of $H$
either to a vertex of
$A := (\{0, s+1\} \times[s]) \cup([s] \times\{0,s+1\})$ or to one of
the at most $\alpha s$ vertices of
$B := [s]^2 \setminus H$.
Observe also that every vertex of $A$ can be adjacent to at most one
vertex of $H$, while
a vertex of $B$ can be adjacent to at most 4 vertices of $H$. Hence, we have
%
\begin{equation}
\label{eq:isosum} e\bigl(H,H^c\bigr) \leq|A| + 4|B| \leq4(1+\alpha)s.
\end{equation}
The claim follows by applying Lemma~\ref{lem:new} [in the form
$v(H)\leq\sqrt{x}\cdot e(H, H^c)/4$] to $H_\beta$.
\end{pf}

%
%
%
%

\begin{lemma}\label{lem:ExtremalkLargest}
Let $\alpha> 0$ and $k \in\eN$ be fixed.
Then the following holds as $s \to\infty$.
If $H \subseteq[s]^2$ is an induced subgraph with $v(H) \geq s^2 -
\alpha s$ vertices,
and $C_1, \ldots, C_k \subseteq H$ denote the $k$ largest components of
$H$ (ties broken arbitrarily) then
\[
v(C_1)+\cdots+v(C_k) \geq\bigl(1-o_s(1)
\bigr) \cdot\lambda_k \cdot s^2,
\]
where $\lambda_k = \lambda_k(\alpha)$ is given by
%
\begin{equation}
\label{eq:lambdadef} \lambda_1 = \biggl(\frac{1}{1+\alpha}
\biggr)^2 \quad\mbox{and}\quad \lambda_{k+1} = \lambda_k
+ \biggl(\frac{1-\lambda_k}{1+\alpha
} \biggr)^2 \qquad\mbox{for } k \geq1.
\end{equation}
\end{lemma}

\begin{pf}
Let us first point out that $0 < \lambda_k < 1$ for all $k$, as can
easily be seen from the definition.
The proof is by induction on $k$. We start with the base case, $k=1$.
Set $\beta:= (1-\eps)\cdot\lambda_1 = (1-\eps) \cdot (\frac
{1}{1+\alpha} )^2$, with $0<\eps<1$ arbitrary but fixed.
By Lemma~\ref{lem:ExtremalMother}, the union of all components of order
at most $\beta s^2$
contains no more than $\sqrt{1-\eps} \cdot s^2$ vertices.
Since the union of \emph{all} components must clearly have $s^2 -
\alpha
s = (1-o_s(1)) s^2$ vertices, there must exist a
component of order $> \beta s^2$.
As $\eps>0$ can be chosen arbitrarily small, it follows that $v(C_1)
\geq(1-o_s(1)) \lambda_1 s^2$, which establishes the base case.

Now suppose that $v(C_1) + \cdots+ v(C_k) = \tilde{\lambda} s^2$ with
$\tilde{\lambda} \geq(1-o_s(1)) \lambda_k$.
If $\tilde{\lambda} > \lambda_{k+1}$, then we are done, so we can
assume this is not the case.
Aiming for a contradiction, suppose that $v(C_{k+1}) < \beta s^2$, where
$\beta= (1-\eps) (\frac{1-\tilde{\lambda}}{1+\alpha}
)^2$ for
some fixed $\eps> 0$.
Lemma~\ref{lem:ExtremalMother} would then give that
\[
v\bigl( H \setminus(C_1\cup\cdots\cup C_k) \bigr) \leq
\sqrt{\beta} (1+\alpha) s^2 = \sqrt{1-\eps} \cdot(1-\tilde{\lambda})
s^2,
\]
which is impossible as we must have
\begin{eqnarray*} \bigl(1-o_s(1)\bigr) s^2
& = & v(C_1\cup\cdots\cup C_k) + v\bigl( H
\setminus(C_1\cup\cdots\cup C_k) \bigr)
\\
& = & \tilde{\lambda} s^2 + v\bigl( H \setminus(C_1\cup
\cdots\cup C_k) \bigr).
\end{eqnarray*}
It follows that $v(C_{k+1}) \geq(1-o_s(1))  (\frac{1-\tilde
{\lambda
}}{1+\alpha} )^2$, so that
\[
v(C_1) + \cdots+ v(C_{k+1}) \geq\bigl(1-o_s(1)
\bigr) \cdot \biggl( \tilde {\lambda } + \biggl(\frac{1-\tilde{\lambda}}{1+\alpha}
\biggr)^2 \biggr).
\]
By differentiating $f(x) := x +  (\frac{1-x}{1+\alpha} )^2$
with respect to $x$, it is easily seen that $f$ is
strictly increasing in $x$ for $x \geq\lambda_1 =  (\frac
{1}{1+\alpha} )^2$.
Since $\tilde{\lambda} \geq(1-o_s(1)) \lambda_k$ by the inductive
hypothesis,
and $\lambda_k > \lambda_{k-1} > \cdots> \lambda_1$, it now follows that
\begin{eqnarray*} v(C_1) + \cdots+
v(C_{k+1}) & \geq& \bigl(1-o_s(1)\bigr) \cdot \biggl(
\lambda_k + \biggl(\frac{1-\lambda
_k}{1+\alpha
} \biggr)^2 \biggr)
\cdot s^2
\\
& = & \bigl(1-o_s(1)\bigr) \cdot\lambda_{k+1} \cdot
s^2,
\end{eqnarray*}
as required.
\end{pf}

\begin{corollary}\label{cor:Extremalk}
Fix $0 < \eps< 1, \alpha> 0$.
Then there exists $k=k(\eps,\alpha)$ such that the following holds for
all large enough $s$.
If $H \subseteq[s]^2$ is an induced subgraph with $v(H) \geq s^2 -
\alpha s$ then
\[
v(C_1) + \cdots+ v(C_k) \geq(1-\eps) s^2,
\]
where $C_i$ denotes the $i$th largest component (ties broken arbitrarily).
Moreover, for all $\eps$ we have $k(\eps,\alpha) = 1$ if $\alpha
\leq
\eps/2$.
\end{corollary}

\begin{pf}
Let the numbers $\lambda_k = \lambda_k(\alpha)$ be as defined
by \eqref
{eq:lambdadef}.

The ``moreover'' part of Corollary~\ref{cor:Extremalk} follows
immediately from Lemma~\ref{lem:ExtremalkLargest} since
$\lambda_1 =  (\frac{1}{1+\alpha} )^2 > 1 - \eps$ if
$\alpha
\leq\eps/2$.

To see that the rest of the corollary also holds, notice that the
numbers $\lambda_k$ form an increasing
sequence that is bounded above by one,
and that the limit of the sequence must be a fixed point of the equation
\[
\lambda= \lambda+ \biggl(\frac{1-\lambda}{1+\alpha} \biggr)^2.
\]
Since the only fixed point is $\lambda=1$, we must have $\lim_{k\to
\infty} \lambda_k = 1$.
Hence, there is a $k=k(\eps,\alpha)$ such that $\lambda_k > 1-\eps$.
The lemma follows.
\end{pf}

\subsubsection{Balls and bins}

In our proof of part (ii) of Theorem~\ref{thm:online},
we will need a minor extension of the lower bound part of Theorem~\ref
{thm:AzarEtAl} that concerns the scenario where there are slightly
fewer balls than bins. This case does not seem to have been treated
explicitly in the literature.
The proof is very similar to the original lower bound proof given
in \cite{AzarEtAl}; we include it here for completeness.

\begin{lemma}\label{lem:ballsbins}
Let $\eps> 0$ be arbitrary, but fixed.
Consider the power of two choices balls and bins process, with $n$ bins and
$m$ rounds where $n/\ln n \leq m \leq n$, and let $M_n$ denote the
maximum load in round $m$.
No matter what strategy the player utilizes, we have
\[
\Pee\bigl( M_n < (1-\eps) \log\log n \bigr) \leq\exp\bigl[ -
n^{1+o(1)} \bigr].
\]
\end{lemma}

\begin{pf}
By obvious monotonicity properties, it suffices to
prove the lemma for the case when $m$, the number of rounds, is exactly
equal to $n/\ln n$. Let us thus assume that $m = n/\ln n$.
Our approach for the proof will be to bound from below, for each $i$,
the number of rounds in which
the player \emph{is forced to} create a bin with $i$ balls in it.

We denote by $N_i(t)$ the number of bins with \emph{at least} $i$ balls
in them after round~$t$. [Note that $N_0(t) = n$ for all $t$.]
Furthermore, we set
\[
k := \bigl\lceil(1-\eps) \log\log n\bigr\rceil,
\]
and write
$\underline{N}(t) = (N_0(t), \ldots, N_k(t))$. Let
\[
\alpha:= \biggl(\frac{m}{4 \cdot k \cdot n} \biggr),\qquad t_i := m
\cdot(i/k) \qquad\mbox{for } i=1,\ldots, k,
\]
and define $c_i$ by
\[
c_i := \alpha^{2^i-1}, \qquad i=1,\ldots, k.
\]
Observe that the $c_i$ satisfy the recurrence relation
%
\begin{equation}
\label{eq:cirecurr} c_{i+1} = \alpha\cdot c_i^2.
\end{equation}
Before proceeding, let us make some further observations about the $c_i$.
Note that $\alpha= (1+\eps) / (4\log\log n \cdot\log n)$, so
certainly $\alpha< 1$ and the $c_i$ are thus decreasing. Moreover, we have
%
\begin{equation}
\label{eq:lnalphaprop} 0 > \ln\alpha\geq- \bigl(1+o(1)\bigr) \ln\ln n.
\end{equation}
Also note that
%
\begin{eqnarray}
\label{eq:ciprop} %
 c_k & = &
\alpha^{2^k - 1}\nonumber
\\
& \geq& \alpha^{ (\log n)^{1-\eps} - 1 }\nonumber
\\
& = & \exp \bigl[ \ln\alpha\cdot\bigl( (\log n)^{1-\eps} - 1\bigr) \bigr]
\\
& \geq& \exp \bigl[ - \bigl(1+o(1)\bigr)\cdot\ln\ln n \cdot(\log
n)^{1-\eps} \bigr]\nonumber
\\
& = & \exp\bigl[ - o(\ln n) \bigr]\nonumber
\\
& = & n^{-o(1)},\nonumber
\end{eqnarray}
where we have used \eqref{eq:lnalphaprop} in the fourth line.

Another key observation is that if in some round both balls fall in
bins with
\emph{exactly} $i-1$ balls in them then a new bin with $i$ balls in it
will be created, regardless of the strategy of the player.
In other words, for all $i, t$ we have
%
\begin{equation}
\label{eq:ballsbinsobs} \Pee\bigl( N_i(t+1) = N_i(t) + 1 |
\underline{N}(t) \bigr) \geq \biggl(\frac{N_{i-1}(t) - N_{i}(t)}{n} \biggr)^2.
\end{equation}
%
%
If for some $t$ we have $N_{1}(t) < c_1 n$, then
the observation \eqref{eq:ballsbinsobs} shows that
\[
\Pee\bigl( N_1(t+1) = N_1(t) + 1 | N_1(t)
< c_1 n \bigr) \geq(1-c_1)^2 \geq
\tfrac{1}2.
\]
It follows that
%
\begin{eqnarray}
\label{eq:N111} %
\Pee\bigl( N_1(t_1)
< c_1 n \bigr) & \leq& \Pee\biggl( \Bi\biggl( t_1,
\frac{1}2 \biggr) < c_1 n \biggr)
\nonumber\\
& \leq& \exp\biggl[ - (t_1/2) \cdot H \biggl(\frac{c_1 n}{t_1/2}
\biggr) \biggr]
\nonumber
\\[-8pt]
\\[-8pt]
\nonumber
& = & \exp\biggl[ - \frac{m}{2 k} \cdot H\biggl(\frac{1}{2} \biggr)
\biggr]
\\
& = & \exp\bigl[ - n^{1-o(1)} \bigr], \nonumber
\end{eqnarray}
where we used the Chernoff bound (Lemma~\ref{lem:chernoff}), together
with the facts
that $t_1 = m/k = n^{1-o(1)}$ and $c_1 = \alpha= m / (4 k n)$.

Let $E_i$ denote the event
\[
E_i := \bigl\{ N_j(t_j) \geq
c_j n \mbox{ for all }1 \leq j \leq i \bigr\}.
\]
Again using the observation \eqref{eq:ballsbinsobs} we see that, for
all $t \geq t_i$:
\[
\Pee\bigl( N_{i+1}(t+1) = N_{i+1}(t) + 1 | E_i,
N_{i+1}(t) < c_{i+1} n \bigr) \geq (c_i-c_{i+1})^2
\geq\tfrac{1}2 c_i^2.
\]
[Here, we use that $N_i(t) \geq N_i(t_i)$ for $t\geq t_i$ by obvious
monotonicity.]
It thus follows that
%
\begin{eqnarray}
\label{eq:Niii} %
 \Pee\bigl( N_{i+1}(t_{i+1})
< c_{i+1} n | E_i \bigr) & \leq& \Pee\bigl( \Bi\bigl(
t_{i+1} - t_i, c_i^2 / 2 \bigr) <
c_{i+1} n \bigr)\nonumber
\\
& \leq& \exp\biggl[ - (m/k) \cdot\bigl(c_i^2/2\bigr)
\cdot H\biggl(\frac{c_{i+1} n}{ (m/k) \cdot
(c_i^2/2)} \biggr) \biggr]
\nonumber
\\[-8pt]
\\[-8pt]
\nonumber
& = & \exp\biggl[ - (m/k) \cdot\bigl(c_i^2/2\bigr) \cdot
H\biggl(\frac{1}2 \biggr) \biggr]
\\
& = & \exp\bigl[ - n^{1-o(1)} \bigr], \nonumber
\end{eqnarray}
where we used that $c_{i+1} = \alpha c_i^2 = (m/ 4kn) \cdot c_i^2$ in
the third line, and
that $c_i \geq c_k = n^{-o(1)}$ by \eqref{eq:ciprop} in the fourth line.
It follows that
\begin{eqnarray*} \Pee\bigl( N_k(t_k) <
c_k n \bigr) & \leq& \Pee\bigl( \mbox{there is some $1 \leq i \leq k$
such that $N_i(t_i) < c_i n$}\bigr)
\\
& = & \Pee\bigl( E_1^c\bigr) + \Pee\bigl(
E_2^c | E_1 \bigr)\Pee( E_1) +
\cdots+ \Pee\bigl( E_k^c | E_{k-1} \bigr)\Pee(
E_{k-1} )
\\
& \leq& \Pee\bigl( E_1^c\bigr) + \Pee\bigl(
E_2^c | E_1 \bigr) + \cdots+ \Pee\bigl(
E_k^c | E_{k-1} \bigr)
\\
& = & \Pee\bigl( N_1(t_1) < c_1 n \bigr)
+ \sum_{i=1}^{k-1} \Pee\bigl(
N_{i+1}(t_{i+1}) < c_{i+1} n | E_i
\bigr)
\\
& = & k \cdot\exp\bigl[ - n^{1-o(1)} \bigr]
\\
& = & \exp\bigl[ - n^{1-o(1)} \bigr],
\end{eqnarray*}
where we used \eqref{eq:N111} and \eqref{eq:Niii} to get the fifth line.
This concludes the proof of the lemma.
\end{pf}

\subsubsection{The two choices coupon collector}
\label{sec:ccc}

In the classical coupon collector problem, each box of some product
contains one of $N$ types of coupons, sampled uniformly at random.
There is a collector who keeps buying new boxes until he has collected
at least one coupon of each type.

For our purposes, it is useful to consider a variant of the coupon
collector problem that we call the \emph{two choices coupon collector}.
Suppose again that there are $N$ types of coupons, but now each box
contains \emph{two} coupons
(whose types are chosen independently and uniformly at random from all
$N$ possible types).
There is a collector (2CCC, henceforth) who only adds a coupon to his
collection if \emph{both} coupons in the box are of types he does not
have in \mbox{his collection} yet. Even in that case, he only adds one of the
two coupons to his collection. 

\begin{lemma}\label{lem:CCC}
For any $s=s(N)$ with $1 \ll s \ll N$, the following is true: w.h.p.
the 2CCC needs to buy at most $2N^2/s$ boxes to collect all but $s$ coupons.
\end{lemma}

\begin{pf}
Let $T$ denote the number of rounds it takes the 2CCC to collect
exactly $N-s$ coupons.
We need to show that
%
\begin{equation}
\label{eq:CCCeq} \Pee\bigl( T > 2N^2/s \bigr) = o(1).
\end{equation}
Observe that $T$ is a sum of independent geometrically distributed
random variables. More precisely,
\[
T = Z_1 + \cdots+ Z_{N-s},
\]
where $Z_i \isd\geom( p_i )$ with $p_i :=  (\frac
{N-i+1}{N} )^2$.
Thus,
\begin{eqnarray*}\Ee T & = & \sum_{i=1}^{N-s}
\frac{1}{p_i} = \sum_{i=1}^{N-s}
\biggl(\frac{N}{N-i+1} \biggr)^2
\\
& = & N^2 \sum_{j=s+1}^N
j^{-2}
\\
& = & N^2 \cdot\bigl(1+o(1)\bigr) \int_{s}^\infty
x^{-2}{\dd}x
\\
& = & \bigl(1+o(1)\bigr) \frac{N^2}{s},
\end{eqnarray*}
where the integral approximation holds due to our assumptions on $s$.
Similarly, we have
%
\begin{eqnarray}
\label{eq:VarTeq} %
\Var T & = & \sum
_{i=1}^{N-s} \frac{1-p_i}{p_i^2} \leq \sum
_{i=1}^{N-s} \frac{1}{p_i^2}\nonumber
\\
& = & \sum_{i=1}^{N-s} \biggl(
\frac{N}{N-i+1} \biggr)^4\nonumber
\\
& = & N^4 \sum_{j=s+1}^N
j^{-4}
\\
& = & N^4 \cdot\bigl(1+o(1)\bigr) \int_{s}^\infty
x^{-4}{\dd}x\nonumber
\\
& = & \bigl(1+o(1)\bigr)\frac{N^4}{3s^3}. \nonumber
\end{eqnarray}
It follows with Chebyshev's inequality that
\[
\Pee\bigl( T > 2N^2/s\bigr) \leq\Pee( T > 1.5 \Ee T ) \leq
\frac{\Var T}{(0.5 \Ee T)^2} = O(1/s)=o(1),
\]
as desired.
\end{pf}


\subsection{Proof of part \textup{(ii)} of
Theorem \texorpdfstring{\protect\ref{thm:online}}{1}}
\label{sec:online.ub.pf}

We now are ready to give the main argument for our upper bound on the
online threshold. Throughout this section, we will consider the boxes
of the dissection $\Dscr_{\rho}$ as defined in \eqref{eq:Drdef}, where
$\rho:= r / \sqrt{5}$.
This time, however, we treat these boxes as vertices of the ordinary
grid graph $[(1/\rho)]^2$.
[Again we assume $(1/\rho)$ is an integer.]
We will denote by $\Occ_\rho(t)$ the subgraph induced by the occupied
boxes after $t$ rounds.

Note that if two points fall anywhere inside two adjacent boxes of
$\Dscr_{\rho}$, by our choice of $\rho$ they are within distance $r$
from each other. Hence, if a set of occupied boxes induces a connected
component of $\Occ_\rho(t)$, then all points inside these boxes belong
to the same connected component of the geometric graph.

\begin{lemma} \label{lem:CCC-consequence}
After $n/2$ rounds, w.h.p. all but $ (\frac{100}{c} )
\cdot
\log\log n \cdot(1/\rho)$ boxes are occupied, no matter how the
player plays.
\end{lemma}

\begin{pf}
Note that to minimize the number of occupied boxes, the player should
play exactly like the two choices coupon collector from Section~\ref{sec:ccc}, with the $N:=(1/\rho)^2=n^{2/3+o(1)}$ boxes playing the role
of the coupons. Let $s:= (\frac{100}{c} ) \cdot\log\log n
\cdot(1/\rho)$, and note that $s=n^{1/3+o(1)}=o(N)$. It follows with
Lemma~\ref{lem:CCC} that w.h.p. after
at most
\[
\frac{2N^2}{s}=\frac{2c}{100 \cdot\rho^3\log\log n} = \frac
{c\cdot
5^{3/2}}{50 \cdot r^3\log\log n} = \frac{5^{3/2}}{50}
\cdot n < n/2
\]
rounds, all but $s$ many boxes are occupied.
\end{pf}

%
%
%
%
%
%
%
%

Let us now fix $0 < \eps< 1/7$, to be determined later.
We dissect the square into \emph{b-blocks} (which stands for ``big
blocks'') consisting of $b\times b$ boxes, where
%
\begin{equation}
\label{eq:bdef} b := \biggl(\frac{\eps\cdot c}{1000} \biggr) \cdot\frac{(1/\rho
)}{\log
\log n}.
\end{equation}
[Again we assume for convenience that $(1/\rho)$ and $(1/b\rho)$ are
integers.]
Thus, there are $z^2$ b-blocks where
%
\begin{equation}
\label{eq:zdef} z = (1/b\rho) = \biggl(\frac{1000}{\eps\cdot c} \biggr) \cdot\log \log
n. 
\end{equation}
It is convenient to consider the b-blocks as vertices of the (ordinary)
grid $[z]^2$.
So a b-block is adjacent to b-blocks that share a side with it, but not
with b-blocks that
share only a corner with it.

\begin{figure}
\centering
\begin{tabular}{@{}cc@{}}

\includegraphics{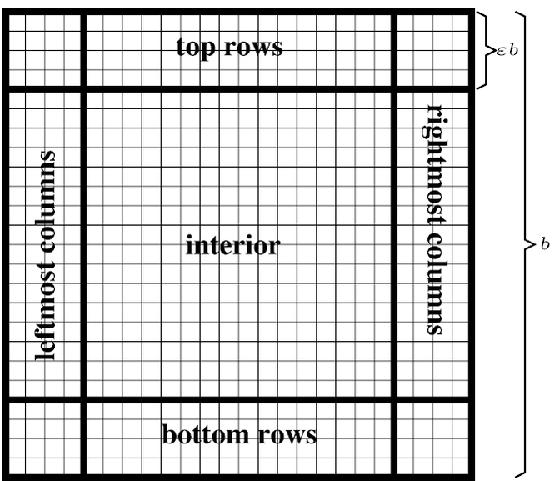}
 & \includegraphics{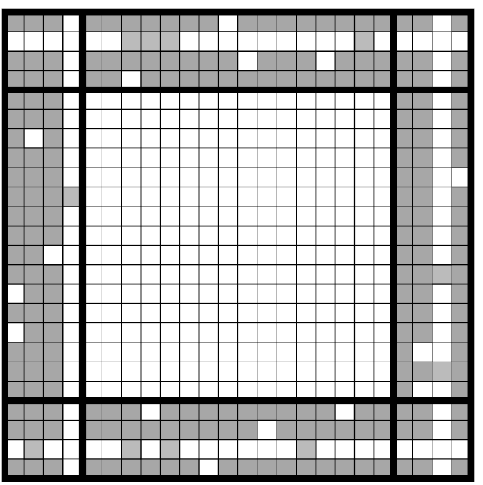}\\
\footnotesize{(a)} & \footnotesize {(b)}\\
\end{tabular}
\caption{Two b-blocks. \textup{(a)} A b-block with its top, bottom, leftmost
and rightmost columns labelled. \textup{(b)} A good b-block. Three of the four top rows
have
only few empty boxes, and similarly for
the bottom rows, left columns and right
columns.}\label{fig:superblocks}
\end{figure}

We shall refer to the top $\eps b$ rows of a b-block $B$ simply as the
\emph{top rows} of $B$.
Similarly, we call the bottom $\eps b$ rows the bottom rows, the
leftmost $\eps b$ columns the leftmost
columns and the rightmost $\eps b$ the rightmost columns.
Those boxes of a b-block $B$ that belong to neither the top or bottom
rows nor
to the leftmost or rightmost columns will be called the \emph{interior}
of $B$.
See Figure~\ref{fig:superblocks}(a) for a depiction.

%
%

%

Let us call a row of a b-block \emph{good} if no more than
$\frac{1}3 \log\log n$ of its boxes are empty, and similarly we call a
column of a b-block
good if no more than $\frac{1}3 \log\log n$ of its boxes are empty.

We call a b-block \emph{good} if
at least three quarters of the top rows are good, at least three
quarters of the bottom rows are good, at least three quarters of the
leftmost columns are good, and
at least three quarters of the rightmost columns are good.
See Figure~\ref{fig:superblocks}(b) for a depiction.

If a b-block is not good, we will call it \emph{bad}.
Let us denote by $\super_\rho(t)$ the subgraph of the grid of b-blocks
(recall that we treat it like the
ordinary $z\times z$ grid) induced by the b-blocks that are good in
round $t$.


\begin{lemma}\label{lem:superblockextr}
W.h.p. in round $n/2$ at most $ (\frac{10^{5}}{\eps c} )
\cdot
z$ b-blocks are bad, no matter what the player does.
\end{lemma}

\begin{pf}
By Lemma~\ref{lem:CCC-consequence}, we can assume that in roun\label{fig:good}d $n/2$
at most
$ (\frac{100}{c} ) \cdot\log\log n\cdot(1/\rho)$ boxes
are empty.
Each bad b-block contains at least
\[
\frac{1}4 \cdot\eps\cdot b \cdot\frac{1}3 \log\log n = \biggl(
\frac{\eps^2 c}{12\scriptsize{\mbox{,}}000} \biggr) (1/\rho),
\]
empty boxes, because
a quarter of either the top rows or the bottom rows or the leftmost columns
or the rightmost columns is not good.
Hence, the number of bad b-blocks cannot be larger than
\[
\frac{ ({100}/{c} ) \cdot\log\log n \cdot(1/\rho)
}{
 ({\eps^2 c}/{12\scriptsize{\mbox{,}}000} ) (1/\rho)
} < \biggl(\frac{10^8}{\eps^2 c^2} \biggr) \log\log n = \biggl(
\frac{10^{5}}{\eps c} \biggr) \cdot z.
\]
\upqed\end{pf}

We shall also need the following consequence of Lemma~\ref
{lem:CCC-consequence}.

\begin{corollary}\label{cor:bcomponents}
W.h.p. in round $n/2$ at least $(1-\eps) (1/\rho)^2$ boxes are
contained in components of $\Occ_{\rho}(n/2)$
of order strictly larger than $b^2$, no matter what the player does.
\end{corollary}

\begin{pf}
By Lemma~\ref{lem:CCC-consequence}, we can assume that in round $n/2$
there are
at most $ (\frac{100}{c} ) \log\log n (1/\rho)$ empty boxes.
Set
\[
s := (1/\rho), \qquad\alpha:= \biggl(\frac{100}{c} \biggr) \cdot\log \log n,\qquad
\beta:= (b/s)^2= \biggl( \frac{\eps\cdot c}{1000 \log\log n } \biggr)^2.
\]
By Lemma~\ref{lem:ExtremalMother}, w.h.p. the number of occupied boxes
in components of $\Occ_{\rho}(n/2)$ of order at most $b^2$ is at most
$\sqrt{\beta} \cdot(1+\alpha) \cdot s^2\leq\sqrt{\beta} \cdot
2\alpha
\cdot s^2=\eps/5 \cdot(1/\rho)^2$. As only an $o(1)$-fraction of the
boxes is empty, the claim follows.
\end{pf}

\begin{figure}

\includegraphics{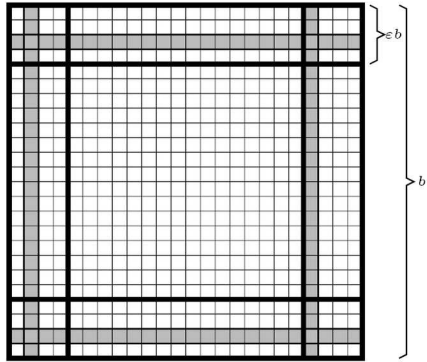}

\caption{A framed b-block.}\label{fig:framed}
\end{figure}

Let us say that a row or a column of a b-block is \emph{full} is it
contains no empty boxes.
We will us say that a b-block is \emph{framed} if
among the top rows there is one that is full, among the bottom rows
there is one
that is full, among the leftmost columns there
is one that is full and among the righmost columns there is one that
is full.
If a b-block $B$ is framed then we refer to the union of the full rows
among the top rows, bottom rows,
leftmost columns and rightmost columns as the \emph{frame} of $B$.
The choice of the name should be clear from the depiction in
Figure~\ref{fig:framed}.

\begin{lemma}\label{lem:framed}
W.h.p. every b-block that is good in round $n/2$ is framed in round~$n$, no matter what the player does.
\end{lemma}

\begin{pf}
We will first compute the probability that a given b-block contains a
full row among the top
$\eps b$ rows in round $n$, given that it was good in round $n/2$.

Let us thus fix a b-block $B$, condition on it being good in round
$n/2$, and
consider what happens to it in the rounds $t > n/2$.
In round $n/2$, at least $M := \frac{3}4 \eps b $ of the top rows of $B$
are good. Let us
fix exactly $M$ of these good rows $r_1, \ldots, r_M$.

We now consider the following balls and bins type process for the
remaining $n/2$ rounds.
In each round $n/2 < t \leq n$, as long as none\vadjust{\goodbreak} of the rows $r_1, \ldots
, r_M$ is full, we have a
list $\underline{e}(t) = (e_1(t),\ldots,e_M(t))$ of empty boxes, where
$e_i(t) \in r_i$.
If at some round $t$ the player plays in some box in our list, then we
replace it with another box
in the same row that is still empty (as long as this is possible; if it
is not possible then evidently
a full row has been created). Otherwise, we keep the list the same.
In other words, if in round $t$ the player picks a point in $e_i(t)$
then we set $e_j(t+1) = e_j(t)$ for all $j\neq i$ and
$e_i(t+1)$ is set to some box in row $i$ that is still empty; and if
the player does not play in any box of $\underline{e}(t)$ then we set
$\underline{e}(t+1) = \underline{e}(t)$.

We want to compute the probability that there is some index $i$ such
that the player is forced to play more than
$\frac{1}3 \log\log n$ times in $e_i(t)$.
Let $R$ denote the number of rounds $n/2 < t \leq n$ in which both points
fall into $e_1(t) \cup\cdots\cup e_M(t)$---so that the player is
forced to play in one of the $e_i$s.
Then
\[
R \isd\Bi\bigl( n/2, \bigl(M \cdot\rho^2\bigr)^2 \bigr).
\]
Observe that
\[
\Ee R = \Theta\bigl( n M^2 \rho^4 \bigr) = \Theta\bigl(
n \cdot\rho^2 / (\log \log n)^2 \bigr) =
n^{{1}/3-o(1)}.
\]
Hence, the Chernoff bound (Lemma~\ref{lem:chernoff}) yields
\[
\Pee( R < \Ee R / 2 ) \leq e^{ - n^{{1}/3-o(1)}}.
\]
Set $N := \frac{1}2\Ee R$.
If we condition on $R > N$, the probability that the player can achieve
a situation where none of $r_1, \ldots, r_M$ is full by round $n$
is upper bounded by the probability that in the player version of the
two choices balls and bins process with $N$ rounds and $M$ bins,
the player can achieve a maximum load of less than $\frac{1}3\log\log n$.
Observe that
\[
M = \tfrac{3}4 \eps b = \Theta\bigl( (1/\rho) / \log\log n \bigr)
= n^{{1}/3-o(1)}.
\]
Thus, we have
\[
\log\log n = \bigl(1+o(1)\bigr) \log\log M.
\]
Also observe that
\begin{eqnarray*} \frac{N}{M} & = & \Theta\bigl( n M
\rho^4 \bigr)
\\
& = & \Theta\bigl( \bigl(n \rho^3\bigr) M \rho\bigr)
\\
& = & \Theta\biggl( \biggl(\frac{c}{\log\log n}\biggr) \cdot\biggl(
\frac{1}{\log\log n} \biggr) \biggr)
\\
& = & \Theta\bigl( (\log\log M)^{-2} \bigr).
\end{eqnarray*}
Hence, by Lemma~\ref{lem:ballsbins} we have
\begin{eqnarray*}&&\Pee( \mbox{none of $r_1,\ldots,
r_M$ is full in round $n$ } | \mbox{ $r_1,\ldots,
r_M$ were good in round $n/2$} )
\\
&&\qquad\leq \Pee\bigl( R < \tfrac{1}2\Ee R \bigr) + \exp\bigl[ -
M^{1-o(1)} \bigr]
\\
&&\qquad\leq e^{ - n^{{1}/3-o(1)} } + e^{-n^{{1}/3-o(1)}}
\\
&&\qquad= e^{ - n^{{1}/3-o(1)} }.
\end{eqnarray*}
The same argument and computations apply to the bottom $\eps b$ rows,
the leftmost $\eps b$ columns and the
rightmost $\eps b$ columns.
Since there are $z^2 = O ( (\log\log n)^2  )$ b-blocks in
total, the union bound gives us
\begin{eqnarray*}&&\Pee(\mbox{there is a b-block which is good in
round $n/2$ and not framed in round $n$} )
\\
&&\qquad\leq z^2 \cdot4 \cdot\exp\bigl[ - n^{{1}/3-o(1)} \bigr]
\\
&&\qquad= o(1),
\end{eqnarray*}
as required.
\end{pf}

\begin{figure}

\includegraphics{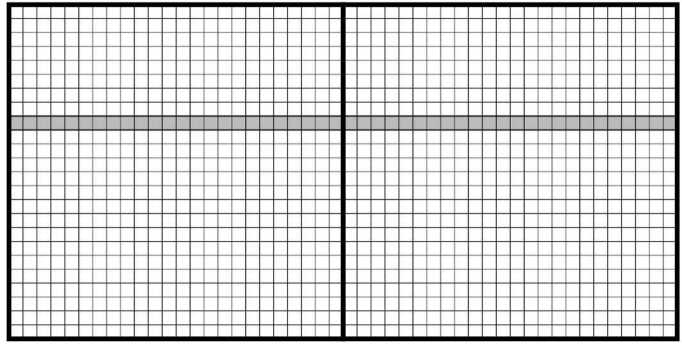}

\caption{Two adjacent, skewered b-blocks.}\label{fig:impaled}
\end{figure}

Let us say that two b-blocks $B_1, B_2$ that share a vertical side are
\emph{skewered} if
there is a row that is full in both $B_1$ and $B_2$.
Similarly, we say that two b-blocks $B_1, B_2$ that share a horizontal
side are skewered if there
is a common column that is full in both. See Figure~\ref{fig:impaled}
for a depiction.

\begin{lemma}\label{lem:impaled}
W.h.p. every two adjacent b-blocks that are both good in round
$n/2$ are skewered in round $n$, no matter what the player does.
\end{lemma}

\begin{pf}
Let $B_1, B_2$ be two adjacent b-blocks (without loss of generality we
can assume
they share a vertical side).
If we condition on both being good in round $n/2$, then there must be
at least
$M = 2\cdot\frac{1}2 \eps b =\eps b$ rows that are good in both $B_1, B_2$.
A row that is good in both has at most $\frac{2}3 \log\log n$ empty boxes.
We can thus follow the same reasoning as in the proof of Lemma~\ref
{lem:framed} and
the same computations with only very minor adaptations to prove that, with
probability $1 - \exp[ - n^{{1}/3-o(1)} ]$ at least one of these rows
will be full in
both $B_1$ and $B_2$.
Since there are only $O(z)=O(\log\log n)$ pairs of adjacent b-blocks, the
union bound again completes the proof.
\end{pf}

Observe that, if two adjacent b-blocks $B_1, B_2$ are both framed and
if they are also skewered, then
their frames will belong to the same component of $\Occ_\rho$.
Hence, Lemmas \ref{lem:framed} and \ref{lem:impaled} together
immediately imply the following.

\begin{corollary}\label{cor:superblocksconn}
W.h.p. the folllowing holds, no matter what the player does.
If $\Ccal$ is a connected component of $\super_\rho(n/2)$, then
every b-block of $B \in\Ccal$ will be framed in round $n$ and the frames
will all belong to the same component of the boxes graph $\Occ_\rho
(n)$. 
\end{corollary}

\begin{lemma}\label{lem:ac}
W.h.p., no matter what the player does, in round $n$ there will be a component
$\Cscr$ of $\Occ_{\rho}(n)$ that contains
at least $v(\Cscr) \geq a(\eps,c) \cdot(1/\rho)^2$ boxes, where
%
\begin{equation}
\label{eq:acdef} a(\eps,c) = \frac{1-7\eps}{k (\eps, {10^{5}}/{(\eps
c)} )},
\end{equation}
with $k(\cdot,\cdot)$ as provided by Corollary~\ref{cor:Extremalk}.
\end{lemma}


\begin{pf}
Let $A_1 \subseteq[0,1]^2$ denote the union of all boxes that
belong to a component of $\Occ_{\rho}(n/2)$ of order larger than $b^2$
(in round $n/2$).
By Corollary~\ref{cor:bcomponents}, w.h.p., we have
\[
\area(A_1) \geq1 - \eps.
\]
Let $A_2$ denote the intersection of $A_1$ with the union of all boxes that
lie in the interior of some good b-block $\in\super_\rho(n/2)$. %
(The reason for these definitions will become clear later.)
Then
\begin{eqnarray*} \area(A_2) & \geq&
\area(A_1) - 4\eps(b\rho)^2 N_{\mathrm{good}} - (b
\rho)^2 N_{\mathrm
{bad}}
\\
& \geq& \area(A_1) - 4\eps- (b\rho)^2 N_{\mathrm{bad}},
\end{eqnarray*}
where $N_{\mathrm{good}}, N_{\mathrm{bad}}$ denote the number of good
respectively  bad b-blocks in
round $n/2$. Trivially, we have $N_{\mathrm{good}} \leq z^2 =(1/b\rho
)^2$. Moreover,
by Lemma~\ref{lem:superblockextr} we have that $N_{\mathrm{bad}} = O( z
) = o( z^{2} )$ w.h.p.
Hence, w.h.p. it holds that
\[
\area( A_2 ) \geq\area( A_1 ) -5\eps\geq1 - 6 \eps.
\]
Recall that, by Lemma~\ref{lem:superblockextr}, $\super_\rho(n/2)$
contains at least $z^2 - \alpha z$ b-blocks, where $\alpha:= 10^{5} /
(\eps c)$.
By Corollary~\ref{cor:Extremalk}, the $k = k(\eps, \alpha)$ largest
components of $\super_\rho(n/2)$ together
cover a fraction of at least $1-\eps$ of the unit square.
Let $A_3$ denote the intersection of $A_2$ with the b-blocks belonging
to the
$k$ largest components of $\super_\rho(n/2)$.
Then
\[
\area(A_3) \geq\area(A_2) - \eps\geq1 - 7\eps.
\]
Let $\Ccal_1, \ldots, \Ccal_k$ denote the $k$ largest components of
$\super_\rho(n/2)$.
Denoting by $B_i$ the union of the b-blocks of $\Ccal_i$ for each $i=1,
\ldots, n$, there is an index $1\leq i \leq k$ such that
\[
\area( A_3 \cap B_i ) \geq\area(A_3) / k
\geq(1-7\eps) / k.
\]
Now recall that $A_3$ is a union of boxes that belong to components of
$\Occ_\rho(n/2)$ of order at least $b^2$, and that all these boxes
belong to the interior of a b-block belonging to $\Ccal_i$.
Observe that, if $B \in\super_\rho(n)$ is a framed b-block, and
$\Cscr
\subseteq\Occ_\rho(n)$ is a
component with more than $b^2$ boxes that intersects the interior of
$B$, then $\Cscr$ also intersects the frame of $B$
(as a b-block consists of exactly $b^2$ boxes).

By Corollary~\ref{cor:superblocksconn}, we can assume that the frames
of the b-blocks in $\Ccal_i$ all belong to the
same component of the box graph.
Hence, all boxes of $A_3 \cap B_i$ belong to the same component of
$\Occ
_\rho(n)$.
Consequently, $\Occ_\rho(n)$ has a component consisting of at least
$ ( (1-7\eps) / k  ) \cdot(1/\rho)^2$
boxes, as required.
\end{pf}

To transfer the result back to the original random geometric graph
setting and conclude the proof, we need the next lemma, which is
similar in spirit to Lemma~\ref{lem:online.lb.area}.

\begin{lemma}\label{lem:online.ub.area}
For every $\eps> 0$, the following holds w.h.p.
Every $A \subseteq[0,1]^2$ that is the union of boxes of the
dissection $\Dscr_\rho$
and with $\area(A) \geq\eps$ contains at least $(1-\eps)\cdot\area
^2(A) \cdot n$ points
in round $n$, no matter what the player does.
\end{lemma}

\begin{pf}
Let $\Ascr$ denote the set of all sets $A\subseteq[0,1]^2$ under
consideration.
Since every $A \in\Ascr$ is a union of boxes, we have $|\Ascr| \leq
2^{(1/\rho)^2} = 2^{n^{2/3+o(1)}}$.

Fix a set $A \in\Ascr$, and let $Z$ the number of rounds in which the
player cannot avoid playing
in $A$ because both points fall inside it.
Then $Z \isd\Bi( n, \area^2(A) )$.
By the Chernoff bound (Lemma~\ref{lem:chernoff}), we have
\[
\Pee\bigl( Z < (1-\eps)\cdot\area^2(A)\cdot n \bigr) \leq \exp
\bigl[ - n \cdot\area^2(A) \cdot H ( 1-\eps ) \bigr] = \exp\bigl[ -
\Omega(n) \bigr].
\]
This holds for every set $A$ under consideration.
Hence, by the union bound we have
\begin{eqnarray*} &&\Pee\bigl(\mbox{there is a set $A \in\Ascr$
that receives less than $(1-\eps) \cdot\area^2(A) \cdot n$ points}
\bigr)
\\
&&\qquad\leq 2^{ n^{{2}/3+o(1)}} \cdot\exp\bigl[ - \Omega(n) \bigr]
\\
&&\qquad= \exp\bigl[ n^{{2}/3+o(1)} - \Omega(n) \bigr]
\\
&&\qquad= o(1),
\end{eqnarray*}
which gives the lemma.
\end{pf}

With Lemmas \ref{lem:CCC}--\ref{lem:online.ub.area} in hand, it is easy
to prove part (ii) of Theorem~\ref{thm:online}.

\textit{Proof of part \textup{(ii)} of
Theorem~\ref{thm:online}}:
For given $c$, set $\eps:=\sqrt{2\cdot10^5/c}>0$ if $c\geq10^8$, and
$\eps=0.01$ otherwise. Let $\widetilde a(c):=a(\eps,c)$ for $a(\eps,c)$
as defined in \eqref{eq:acdef}. Note that in both cases $0<\widetilde
{a}(c)<1$. Further, by the ``moreover'' part of Corollary~\ref
{cor:Extremalk}, for $c\geq10^8$ we have $\widetilde{a}(c)=1-7\eps
=1-O(1/\sqrt{c})$. Hence, we have $\widetilde{a}(c)\to1$ as $c\to
\infty$.

By Lemma~\ref{lem:ac}, w.h.p. the boxes graph $\Occ_\rho(n)$ will have
a connected component of area at least
$\widetilde{a}(c)$. By Lemma~\ref{lem:online.ub.area}, this will give
us a component of order at least $(1-o(1)) \cdot\widetilde a^2(c)
\cdot n$ in the resulting geometric graph. Hence, the claim follows
for, say, $g(c) := \widetilde a^3(c)$. Note that $g(c)\to1$ as $c\to
\infty$.


\subsection{Proof of part \textup{(ii)} of
Theorem \texorpdfstring{\protect\ref{thm:offline}}{2}}
\label{sec:offline.ub.pf}


As in the previous proof, we divide the unit square into boxes of
side-length $\rho:=r/\sqrt{5}$. We set $s:=1/\rho$ as before
and assume, also as before, that $s$ is an integer.
Again we denote by $\Occ_\rho$ the subgraph of the $s \times s$ grid
$[s]^2$ induced by the occupied boxes.

For given $c>0$, define $a=a(c)$ as the solution of
%
\begin{equation}
\label{eq:def-eps} \frac{480\sqrt{a}}{(1-\sqrt{a})^2} =c.
\end{equation}
Note that $0<a(c)<1$ with $a(c)\to1$ as $c\to\infty$.

We will show the following.

\begin{claim} \label{claim:grid-counting}
W.h.p. every possible choice of points is such that $\Occ_\rho$ has a
component with more than $a(c) \cdot s^2$ vertices.
\end{claim}

Observing that Lemma~\ref{lem:online.ub.area} carries over to the
offline setting, Claim~\ref{claim:grid-counting} implies part (ii) of Theorem~\ref{thm:offline} as in the argument just
given for the online case [for, say, $g(c):=a^2(c)$].

It therefore remains to prove Claim~\ref{claim:grid-counting}. To do
so, we proceed by combinatorial counting in the grid $[s]^2$. Let $\cX$
denote the family of all subsets $X \subseteq[s]^2$ for which all
components of the graph $[s]^2\setminus X$ are of order at most $a s^2$.
Note that a choice of points for which $\Occ_\rho$ has only components
of order at most $a s^2$ exists if and only if there is a set $X\in\cX$
that can be completely avoided by the player; that is, if and only if
there is a set $X\in\cX$ such that in each of the $n$ pairs, at most
one point falls into one of the boxes of $X$.

Naively speaking, we would wish to show that the expected number of
such sets $X$ is $o(1)$. Then Claim~\ref{eq:def-eps} would follow with
Markov's inequality. Unfortunately, the number of sets $X\in\cX$ is too
large for this. We therefore refine our basic idea by defining a more
manageable family $\mathcal{X}^*$ of subsets of $[s]^2$ with the
crucial property that each $X\in\cX$ has a subset $X^* \subseteq X$
with $X^* \in\cX^*$.
For this family $\mathcal{X}^*$, we will indeed be able to show that
the expected number of sets $X^*\in\cX^*$ that can be avoided in the
sense discussed above is $o(1)$. Once this is established, it follows
with Markov's inequality that w.h.p. no set from $\cX^*$ can be
avoided, which in turn implies that also no set from $\cX$ can be
avoided (recall that each set $X\in\cX$ contains a subset $X^* \in
\cX^*$).
To avoid confusion, let us point out explicitly that $\cX^*$ will
\emph
{not} be a subfamily of $\cX$.

In the following, we proceed with the construction of $\cX^*$. For a
given set $X\in\mathcal{X,}$ we denote by $\cC(X)$ the set of all
components of $[s]^2\setminus X$.
Set
%
\begin{equation}
\label{eq:def-delta} \delta:=\frac{1-\sqrt{a}}{4}.
\end{equation}
%
For $X\in\cX$ given, let
%
\begin{equation}
\label{eq:def-s} t\bigl(|X|\bigr):=\frac{\delta^2 s^4}{(|X|+s)^2},
\end{equation}
and let $k=k(X)$ denote the number of components in $\cC(X)$ that are
of size at least $t(|X|)$. We shall refer to these components as \emph
{large} components, and denote them by $C_1, \ldots, C_k$. We call the
remaining components \emph{small}. Note that the notions of large and
small components are not absolute, but depend on the size of the set
$X$ considered.

For the following definitions, it is convenient to go back to a
geometric viewpoint of the $s\times s$ grid.
Each component $C$ of the graph $[s]^2\setminus X$ corresponds to a
connected subset of the unit square [with area $v(C)\cdot\rho^2$], and
has a geometrical boundary $\partial C$ that is the union of one or
several closed (rectilinear) walks in the unit square. Note that the
length of this geometrical boundary is $e(C,C^c)\cdot\rho$, where
$C^c$ denotes the complement of $C$ in $\mathbb{Z}^2$ (not in $[s]^2$).
For brevity we write, with slight abuse of notation, $|\partial C|$ for
$e(C,C^c)$, and $|C|$ for $v(C)$ in the following.


For $i=1, \ldots, k$, let $C_i'$ denote the maximal superset of $C_i$
whose geometrical boundary $\partial C_i'$ is contained in $\partial
C_i$. Note that $\partial C_i'$ is a single closed walk in the unit
square. (Informally speaking, $C_i'$ is obtained from $C_i$ by ``filling
the holes'' in $C_i$.) Note that $C_1',\ldots, C_k'$ are not necessarily
pairwise disjoint (think, e.g., of $C_1, \ldots, C_k$ as concentric rings).

Going back to the combinatorial viewpoint, it is not hard to see that
the neighbourhood of $C_i'$ is contained in the neighbourhood of
$C_i$ for each $i=1, \ldots, k$.
Let $X'=X'(X)$ denote the union of the neighbourhoods of $C_1',\ldots,
C_k'$, and note that $X'\seq X$.

For $X\in\cX$, we now define
\[
X^*=X^*(X):=\cases{ X, &\quad $\mbox{if }|X|\geq s^{1.01}$, \vspace *{2pt}
\cr
X'(X), & \quad $\mbox{otherwise}$.}
\]

Note that $X^*(X)\seq X$ in both cases---as explained above, this is
crucial for our argument.
Finally, we define $\cX^*$ to be the family of all sets $X^*\seq[s]^2$
that can arise in this way from some set $X\in\cX$.

By our explanations above, it remains to show the following.

\begin{claim} \label{clm:expectation}
The expected number of sets $X^*\in\cX^*$ that contain no two points
from the same random point pair (i.e., the expected number of sets
$X^*\in\cX^*$ that can be avoided by the player) is $o(1)$.
\end{claim}

Let $\cX^*_m$ denote the family of all sets $X^*\in\cX^*$ of size
exactly $m$. We will bound the number of sets in $\cX^*_m$ by
combinatorial counting. We begin by showing that $\cX^*_m$ is in fact
empty for values of $m$ smaller than
%
\begin{equation}
\label{eq:def-kmin} \kmin:=\frac{(1-\sqrt{a})s}{2\sqrt{a}}.
\end{equation}

\begin{lemma} \label{lem:x-prime-large}
For $s$ large and $X\in\cX$ with $|X|< s^{1.01}$ the set $X'=X'(X)$
satisfies $|X'|\geq\kmin$. Consequently, for $m<\kmin$ we have $\cX
^*_m=\varnothing$.
\end{lemma}

\begin{pf}
Let $X$ as in the lemma be given, and let $H_{\mathrm{small}}$ denote
the union of the small components in $\cC(X)$ [recall the definitions
after \eqref{eq:def-s}].
Applying Lemma~\ref{lem:ExtremalMother} with $\beta
=(t(|X|)/s)^2=\frac
{\delta^2 s^2}{(|X|+s)^2}$ and $\alpha=|X|/s$ gives that
$v(H_{\mathrm
{small}})\leq\delta s^2$. The remaining occupied boxes must be in the
$k=k(X)$ large components. Consequently, we have
%
\begin{equation}
\label{eq:sum-large-components} \sum_{i=1}^k
|C_i|\geq(1-\delta)s^2-|X| \geq(1-2\delta)s^2,
\end{equation}
where the last inequality follows from $|X|\leq s^{1.01}=o(s^2)$.

The next argument is similar in spirit to the proofs of Lemmas \ref
{lem:new} and \ref{lem:ExtremalMother}; however, we have to deal with
the subtlety that we want to bound $\sum_{i=1}^k |\partial C'_i| =
\sum_{i=1}^k e(C'_i, (C'_i)^c)$ from below but we only have an upper bound
on $|C_i|=v(C_i)$ [not $|C'_i|=v(C'_i)$] for all $i$.

Note first that
%
\begin{equation}
\label{eq:sum-boundaries-prime} \sum_{i=1}^k \bigl|\partial
C_i'\bigr| \leq4\bigl|X'\bigr| + 4s,
\end{equation}
where the inequality follows from the observation that each vertex of
$X'$ contributes at most $4$ to the sum, and the boundary of the unit
square contributes at most $4s$ in total.

On the other hand, by Lemma~\ref{lem:isoperimetric} the total
circumference of $C'_1, \ldots, C'_k$ satisfies
\begin{eqnarray*}
\sum_{i=1}^k \bigl|\partial
C'_i\bigr| &\geq&\sum_{i=1}^k
4\sqrt{\bigl|C'_i\bigr|} \geq4\sum
_{i=1}^k \sqrt{|C_i|}\geq 4\sum
_{i=1}^k \frac{|C_i|}{\sqrt{a s^2}}
\\
& \geBy{eq:sum-large-components}& \frac{4(1-2\delta)s}{\sqrt{a}},
\end{eqnarray*}
where in the third inequality we used that $|C_i|\leq as^2$ (because
$X\in\cX$).

Together with \eqref{eq:sum-boundaries-prime}, it follows that
\[
\bigl|X'\bigr|\geq\frac{(1-2\delta)s}{\sqrt{a}}-s \eqBy{eq:def-delta}
\frac
{(1-\sqrt{a})s}{2\sqrt{a}}=\kmin.
\]
\upqed\end{pf}

Next, we bound the size of $\cX^*_m$ for the intermediate values of $m$.

\begin{lemma} \label{lem:middle-range}
For $s$ large and $\kmin\leq m < s^{1.01}$, we have $|\cX^*_m|\leq\break 
e^{10m/(1-\sqrt{a})}$.
\end{lemma}

\begin{pf}
To specify a set $X^*\in\cX^*_m$ for $m$ as in the lemma, it suffices
to specify the boundaries of $C_1',\ldots, C_k'$ (i.e., the \emph
{outer} boundaries of $C_1,\ldots, C_k$). For a given such component
$C'_i$, we encode its geometric boundary $\partial C_i'$ by specifying,
say, the leftmost point of the topmost horizontal line intersecting
with $\partial C'_i$ as a starting point, and by specifying the
direction of each of the $|\partial C_i'|=\ell_i$ steps along the
boundary (say in clockwise direction). There are at most $s^2\cdot
3^{\ell_i}$ ways of specifying a boundary $\partial C'_i$ (and thus a
component $C_i'$) in this way.

For each set $X^*\in\cX^*_m$ with $m$ as in the lemma there exists, by
definition, a set $X\in\cX$ with $X'(X)=X^*$ and $|X|< s^{1.01}$. Thus,
the number $k=k(X)$ of large components can be bounded as
\[
k \leq\frac{s^2}{t(|X|)}\eqBy{eq:def-s}\delta^{-2} \biggl(
\frac
{|X|}{s}+1 \biggr)^2 \leq\delta^{-2}
\bigl(s^{0.01}+1\bigr)^2=:x.
\]
Note that $x= O(s^{0.02})$.
Recall also from \eqref{eq:sum-boundaries-prime} that $\sum_{i=1}^k
\ell_i \leq4m+4s$.

It follows that for $\kmin\leq m < s^{1.01}$ we have
\[
\bigl|\cX^*_m\bigr|\leq\sum_{k=1}^{x}
\mathop{\sum_{\ell_1, \ldots, \ell_k\dvtx }}_{
\ell_1+\cdots+ \ell_k\leq4m+4s} \prod
_{i=1}^k \bigl(s^2 3^{\ell_i}
\bigr) \leq{x} (4m+4s)^x \bigl(s^2\bigr)^{x}
3^{{4m+4s}}.
\]

Observe that, for $s$ large enough and $\kmin\leq m <s^{1.01}$, the
factor $3^{{4m+4s}}$
is much larger than $x$, $(4m+4s)^x$ and $(s^2)^{x}$ (which are all
$\exp[O(s^{0.02}\log s)]$).
It follows that, for $s$ large enough:
\[
\bigl|\cX^*_m\bigr|\leq(3.01)^{4(m+s)}\leq e^{5(m+s)}\leq
e^{10m/(1-\sqrt{a})},
\]
where in the last step we used that
\[
s\leBy{eq:def-kmin} \frac{2\sqrt{a}}{1-\sqrt{a}}\cdot\kmin\leq \frac
{2\sqrt{a}}{1-\sqrt{a}} \cdot m
\]
and consequently
\[
m+s\leq\frac{1+\sqrt{a}}{1-\sqrt{a}}\cdot m\leq\frac{2m}{1-\sqrt{a}}.
\]
\upqed\end{pf}

With the preceding lemmas in hand, Claim~\ref{clm:expectation} follows
with a routine calculation.

\begin{pf*}{Proof of Claim~\ref{clm:expectation}}
Using Lemmas \ref
{lem:x-prime-large} and \ref{lem:middle-range}, and using the trivial bound
\[
\bigl|\cX^*_m\bigr|\leq\pmatrix{s^2
\cr
m}\leq
s^{2m}=e^{2m\log s}
\]
for $m\geq s^{1.01}$, we obtain that the expected number of sets
$X^*\in
\cX^*$ that contain no two points from the same random point pair is
%
\begin{eqnarray}
\label{eq:expectation} &&\sum_{m}\bigl|\cX^*_m\bigr|\cdot
\biggl(1- \biggl(\frac{m}{s^2} \biggr)^2 \biggr)^{n}\nonumber\\
&&\qquad
\leq\sum_{\kmin\leq m < s^{1.01}} e^{10m/(1-\sqrt
{a})-m^2n/s^4} + \sum
_{m\geq s^{1.01}} e^{2m\log s-m^2n/s^4}
\\
&&\qquad\leq\sum_{m\geq\kmin} \bigl(e^{10/(1-\sqrt{a})-mn/s^4}
\bigr)^m + \sum_{m\geq s^{1.01}}
\bigl(e^{2\log s-mn/s^4} \bigr)^m.\nonumber
\end{eqnarray}
By our choice of constants, we have
%
\begin{equation}
\label{eq:nLcubed} n/s^3=n\rho^3 = nr^3
\cdot5^{-3/2}\geq c/12
\end{equation}
and consequently the exponents of the terms in parentheses are
uniformly bounded by
\[
\frac{10}{1-\sqrt{a}}-\kmin\cdot n/s^4 \leByM{\scriptsize{\eqref{eq:def-kmin},
\eqref{eq:nLcubed}}} \frac{10}{1-\sqrt{a}} - \frac{c(1-\sqrt
{a})}{24\sqrt
{a}} \eqBy{eq:def-eps} -
\frac{10}{1-\sqrt{a}}<0
\]
and
\[
2\log s - s^{1.01} \cdot n/s^4 \leBy{eq:nLcubed} 2\log s -
c/12\cdot s^{0.01}= -\omega(1),
\]
respectively.
It follows that the right-hand side of \eqref{eq:expectation} is
$o(1)$. 
\end{pf*}

As explained above, Claim~\ref{clm:expectation} implies Claim~\ref
{claim:grid-counting}, which in turn
implies part (ii) of Theorem~\ref{thm:offline}.

\section{Concluding remarks}
\label{sec:concl}

In this paper, we have shown that in the power of choices version of
the random geometric
graph, the onset of a giant component can be delayed until the
average degree is of order $n^{1/3}(\log\log n)^{2/3}$. This is an
improvement by a power of $n$ over the standard random
geometric graph, where a giant appears as soon as the average degree
exceeds a certain constant. As pointed out in the \hyperref[sec1]{Introduction}, this
behaviour is in stark contrast to what happens in the (vertex)
Achlioptas process, where the power of choices only yields a constant
factor improvement.

We have also shown that in the \emph{offline} version of our process a
giant can be delayed
just a little longer, until the average degree is of order $n^{1/3}$.

We offer the following two natural conjectures.

\begin{conjecture}
There is a function $f \dvtx (0,\infty) \to(0,1)$ such that the following holds.
Consider the \emph{online} power of choices geometric graph process,
where $r = \sqrt[3]{\frac{c}{n\log\log n}}$.
Assuming optimal play, the largest component will have size
$(1+o(1))\cdot f(c) \cdot n$ w.h.p.
\end{conjecture}

\begin{conjecture}
There is a function $f \dvtx (0,\infty) \to(0,1)$ such that the following holds.
Consider the \emph{offline} power of choices geometric graph setting,
where $r = \sqrt[3]{\frac{c}{n}}$.
Assuming optimal play, the largest component will have size
$(1+o(1))\cdot f(c) \cdot n$ w.h.p.
\end{conjecture}

Many steps in our proofs have been rather crude and we have made no
attempt to optimize
the expressions for $f(c), g(c)$ in Theorems \ref{thm:online} and \ref
{thm:offline}.
The main reason for this is that we believe that it will not be possible
to prove the above two conjectures without significant new ideas.

\textit{The largest component just before the threshold.}
Our proofs also give some insight into the behaviour before the threshold.
For the following discussion, we let $r_0=(n\log\log n)^{-1/3}$ for
the online case, and $r_0=n^{-1/3}$ for the offline case.
Furthermore, we assume that $r$ is asymptotically smaller than $r_0$
but only slightly so (say $n^{-1/3-0.01}\ll r \ll r_0$).

The strategies described in Sections \ref{sec:online.lb.pf} and
\ref{sec:offline.lb.pf} guarantee that w.h.p. the largest
component is of order $O((r/r_0)^6\cdot n)$ vertices in both settings.
To see this, note that in Lemma~\ref{lem:moat} the constant $a(K)$ can
be improved to $a(K)=\Theta(1/K^2)$. Thus, $a(K)=\Theta(c^2)$ as
$c=(r/r_0)^3\to0$, which translates to the claimed bound by (a~slightly adapted version of) Lemma~\ref{lem:online.lb.area}.

On the other hand, the upper bound proof given in Section~\ref{sec:offline.ub.pf} for the offline setting shows that w.h.p. the
player will be forced to create a component with $\Omega
((r/r_0)^{12}\cdot n)$ vertices: We have $a(c)=\Theta(c^2)$ as
$c=(r/r_0)^3\to0$ in \eqref{eq:def-eps}, and the resulting factor of
$(r/r_0)^6$ is squared when (a~slightly adapted version of) Lemma~\ref
{lem:online.ub.area} is applied.

Similarly, the upper bound proof given in Section~\ref{sec:online.ub.pf} for the online setting shows that w.h.p. the player
will be forced to create a component with at least $\Theta
((r/r_0)^{12}\cdot n)$ vertices: For $\eps$ fixed and $\alpha\to
\infty
$, we have $k(\eps,\alpha)=O(\alpha^2)$ in Corollary~\ref
{cor:Extremalk}, as $\lambda_i=O(i/\alpha^2 )$. It follows that for
$\eps=0.01$ fixed, $a(c,\eps)$ in Lemma~\ref{lem:ac} is $\Omega(c^2)$
as $c=(r/r_0)^3\to0$. As for the online setting, the resulting factor
of $(r/r_0)^6$ is squared when Lemma~\ref{lem:online.ub.area} is applied.

To summarize, in both the online and the offline power of choices
setting, the size of the largest component in optimal play is between
$\Theta((r/r_0)^{12}\cdot n)$ and $\Theta((r/r_0)^6\cdot n)$, where
$r_0$ denotes the respective threshold. Note that this behaviour is
again very different from what happens in the standard geometric and
Erd\H{o}s--R\'{e}nyi random graphs, where the size of the largest
component jumps from $\Theta(\log n)$ to $\Theta(n)$ at the threshold.
(Such a jump is also observed in the Achlioptas process when played
with certain natural, but most likely not optimal player strategies;
see, e.g., \cite{BirthControl}.)

It would be interesting to close the gap between our bounds for the
moderately subcritical regime.

\begin{question}
Both for the online and the offline setting, what is the order of the
largest component in optimal play (w.h.p.) when $r$ is slightly below
the respective threshold?
\end{question}

(Here, we mean by ``slightly below'' that $r_0 \cdot n^{-\eps} \ll r
\ll
r_0$ for every fixed $\eps>0$.)

\textit{More choices.}
Let us now sketch how our results generalize to the scenario with an
arbitrary fixed number $d\geq2$ of choices per step. As stated in the
\hyperref[sec1]{Introduction}, the resulting thresholds then are $n^{-1/(d+1)} (\log
\log
n)^{-(d-1)/(d+1)}$ for the online setting, and
$n^{-1/(d+1)}$ for the offline case. This can be shown with only minor
modifications to the proofs we gave for $d=2$.

To give some intuition for these formulas, let us point out the
following: In both scenarios, the threshold corresponds to the point
where the \emph{number of points that are forced to be in the barrier}
(as defined in our lower bound proofs) equals the \emph{number of boxes
of the barrier} in order of magnitude. In the online scenario, the
barrier has an area of $\Theta(r\log\log n)$, and thus the (expected)
number of points that we need to choose in the barrier is of order $n
r^{d} (\log\log n)^d$. On the other hand, the number of boxes in the
barrier is of order $r^{-1} \log\log n$. It is not hard to see that
these terms are equal for $r$ as stated. The threshold for the offline
case can be motivated with a very similar calculation.

\textit{Creating a giant.}
Another interesting related question is what happens if the
player attempts to speed up the onset of a giant component instead of
delaying it.
For this setup, one can quite easily derive the following result.

\begin{theorem}\label{thm:fertility}
Suppose that $r = \sqrt{\lambda/ n}$ for some constant $\lambda> 0$.
Then the following holds, where $\lambda_{\mathrm{crit}}$ is the
critical constant
for the emergence of a giant component in the ordinary random geometric graph:
\begin{longlist}[(ii)]
\item[(i)] If $\lambda\leq\lambda_{\mathrm
{crit}}/2$ then
w.h.p. the largest component of the graph will be $o(n)$, no matter
what the player does.
\item[(ii)] If $\lambda> \lambda_{\mathrm{crit}} /
2$ then
the player has a strategy that will result in a component of order
$\Omega(n)$, w.h.p.
\end{longlist}
\end{theorem}

To see that part (i) holds, we just need to note that
even if we allow
the player to keep \emph{both} points in each round he will just have
a subcritical or critical random geometric graph.
The proof of part (ii) is only slightly more involved.
A sketch of the argument is as follows:
The player fixes a small square $A$ of area $\eps=\eps(\lambda)$
inside the unit square, and
he always selects a point inside $A$ if he can (if both fall in $A$ he
chooses randomly).
If $\eps>0$ was chosen sufficiently small, then the graph induced by
the points in $A$ will
be a \emph{supercritical} random geometric graph, containing
a linear proportion of all vertices that fall in $A$, and hence also a
linear proportion
of all $n$ vertices.
For completeness, we spell out this argument in more detail in
the Appendix \ref{sec:fertilityapp}.

Our strategy for $\lambda> \lambda_{\mathrm{crit}}/2$ case is rather
simple, and in a sense it might be suboptimal.
While it does deliver a component of linear size, a more sophisticated
strategy might achieve an even larger largest component.

\begin{question}\label{qu:qufert}
If $\lambda> \lambda_{\mathrm{crit}}$ is fixed and $r = \sqrt
{\lambda
/n}$, what is the order
of the largest component the player can (w.h.p.) achieve?
\end{question}

\begin{appendix}\label{app}
\section{An upper bound for the vertex Achlioptas process}
\label{sec:vertex-achlioptas}

In this section, we show that, as claimed at the end of Section~\ref{sec:Background}, in the \emph{vertex} Achlioptas process the player
is also
forced to create a linear-sized component as soon as the average degree
of the underlying random graph exceeds a certain constant.

We will use the following lemma, which is a straightforward
generalization of Lemma~2 in \cite{AvoidingUB}. We denote by $G(n,m)$
a (``Erd\H{o}s--R\'{e}nyi'') random graph sampled uniformly from all
graphs on $n$ vertices and $m$ edges. For a graph $G$ and a set
$S\subseteq V(G)$, we denote the graph induced by $S$ in $G$ by $G[S]$.

\begin{lemma}[(\cite{AvoidingUB})]\label{lem:density_linear}
Let $c > 0$. For every $\eps> 0$ there exists $\delta= \delta(c,
\varepsilon) > 0$ such that \aas the random graph $G := G(n, cn)$ has
the property that for every $S \subseteq V(G)$ for which $G[S]$
contains more than $(1 + \varepsilon)|S|$ edges we have $|S| \geq
\delta n$.
\end{lemma}

We are now able to deduce the following.

\begin{theorem} There is a constant $c>0$ such that if $m \geq cn$ then
w.h.p. a component
of linear size will be formed in the vertex Achlioptas process, no
matter what the player does.
\end{theorem}

\begin{pf}
We will show that for $c$ large enough and $m:=cn$, w.h.p. $G(n,m)$ is
such that every set of $n/2$ vertices induces a graph that contains a
linear-sized component. Clearly, this then proves the claim. (In fact,
our argument gives an upper bound for the \emph{offline} problem
corresponding to the vertex Achlioptas process.)

Note that the expected number of edges in a fixed set of $n/2$ vertices is
\[
m\cdot\frac{{n/2\choose2}}{{n\choose2}}=\bigl(1+o(1)\bigr)m/4.
\]

By a Chernoff-type bound (Theorem~2.10 of \cite{randomgraphs}), the
probability that this number of edges is less than $m/8$ is $e^{-\Omega
(m)}$. A union bound over all (trivially at most $2^n$) sets of $n/2$
vertices thus yields that with probability $1-2^n e^{-\Omega(m)}$, each
such set contains at least $m/8$ edges. Clearly, for $c$ chosen large
enough the last probability is $1-o(1)$.

Note that the ratio of edges to vertices in each such set is
$(m/8)/(n/2)=m/(4n)=c/4$, which we can ensure to be at least 2, say, by
choosing $c\geq8$. Moreover, by averaging, also at least one of the
\emph{components} of the graph induced by such a set has a ratio of
edges to vertices of at least $2$. By Lemma~\ref{lem:density_linear},
w.h.p. each such subgraph of $G(n,m)$ is of order at least $\delta(c,1)n$.

To summarize, w.h.p. $G(n,m)$ is such that each set of $n/2$ vertices
has a ratio of edges to vertices of at least $2$, and as a consequence
of this induces a graph which contains a linear-sized component.
\end{pf}

\section{Proof of the second part of Theorem~
\texorpdfstring{\lowercase{\protect\ref{thm:fertility}}}{37}}
\label{sec:fertilityapp}

In this section, we fill in the details of the proof sketch provided
just after the statement of Theorem~\ref{thm:fertility}.

\textit{Proof of part \textup{(ii)}
of Theorem~\ref{thm:fertility}}:
We take $\eps=\eps(\lambda)$ sufficiently small, to be made more
precise later,
and we let $A \subseteq[0,1]^2$ be a square with $\area(A) = \eps$.

In every round, the player will always picks a point in $A$ if he can.
If it happens that both points fall in $A$ the he chooses randomly.
Observe that the probability that, in a given round, the player is able
to select
a point of $A$ equals $1 - (1-\eps)^2 = 2\eps-\eps^2$.

Let $R$ denote the number of rounds when he succeeded to pick a point
of $A$.
Clearly, $R \isd\Bi( n, 2\eps-\eps^2)$.
By the Chernoff bound (Lemma~\ref{lem:chernoff}), we have that
\[
\Pee\bigl( R < (1-\eps)\cdot\Ee R \bigr) \leq\exp\bigl[ - \Omega( n ) \bigr] =
o(1).
\]
Let $\Gtil$ denote the subgraph of the player's graph induced by the
points in $A$.
Observe that we can rescale $A$ by a factor of $1/\sqrt{\eps}$ and
translate it to
map it to the unit square $[0,1]^2$.
Thus, by stopping the process the instant $n' := (1-\eps)\Ee R =
(1-\eps
)\cdot(2\eps-\eps^2)\cdot n$
points have been selected inside $A$, we see that (w.h.p.) $\Gtil$
contains a
copy of the ordinary random geometric graph with parameters $n'$ and
$r' := r/\sqrt{\eps}$.

Let us now observe that we can rewrite
$r'$ as
\[
r' = \frac{r}{\sqrt{\eps}} = \sqrt{\frac{\lambda}{\eps n}} = \sqrt{
\frac{2 \cdot(1-\eps)\cdot(1-\eps/2) \cdot\lambda}{n'}} =: \sqrt{\frac{\lambda'}{n'}}.
\]
As $\lambda> \lambda_{\mathrm{crit}}/2$, we can choose $\eps>0$ small
enough for $\lambda' > \lambda_{\mathrm{crit}}$ to hold.
Hence, in that case $\Gtil$ will (w.h.p.) contain a component spanning
$\Omega(n') = \Omega(n)$ points.
\end{appendix}

%

%





\printaddresses
\end{document}